\documentclass[11pt,reqno]{amsart}

\usepackage{amsmath,amsthm,amscd,amssymb,amsfonts, amsbsy}
\usepackage{latexsym}
\usepackage{txfonts}
\usepackage{exscale}
\usepackage{mathrsfs}
\usepackage{esint} 
\usepackage{enumitem} 
\usepackage{color}

\usepackage{relsize}
\usepackage[normalem]{ulem}

\usepackage{graphicx} 
\def\Yint#1{\mathchoice
    {\YYint\displaystyle\textstyle{#1}}%
    {\YYint\textstyle\scriptstyle{#1}}%
    {\YYint\scriptstyle\scriptscriptstyle{#1}}%
    {\YYint\scriptscriptstyle\scriptscriptstyle{#1}}%
      \!\iint}
\def\YYint#1#2#3{{\setbox0=\hbox{$#1{#2#3}{\iint}$}
    \vcenter{\hbox{$#2#3$}}\kern-.51\wd0}}
\def\longdash{{-}\mkern-3.5mu{-}} 
\def\tiltlongdash{\rotatebox[origin=c]{15}{$\longdash$}}
\def\fiint{\Yint\tiltlongdash}

\usepackage[
  hmarginratio={1:1},     
  vmarginratio={1:1},     
  textwidth=400pt,			  
	textheight=580pt,        
  heightrounded,          
]{geometry}

\usepackage[utf8]{inputenc}

\usepackage[
	colorlinks,
	citecolor=blue,
	linktoc=all,
	pagebackref,
	hypertexnames=false
]{hyperref}

\usepackage{appendix}

\definecolor{br}{rgb}{1, 0.4,0}

\allowdisplaybreaks

\numberwithin{equation}{section}

\theoremstyle{plain}
\newtheorem{theorem}[equation]{Theorem}
\newtheorem{prop}[equation]{Proposition}
\newtheorem{corollary}[equation]{Corollary}
\newtheorem{lemma}[equation]{Lemma}

\theoremstyle{definition}
\newtheorem{defn}[equation]{Definition}

\theoremstyle{remark}
\newtheorem{remark}[equation]{Remark}

\numberwithin{equation}{section}

\newcommand{\RR}{{\mathbb{R}}}

\newcommand{\HH}{\mathfrak{H}}

\newcommand{\E}{\mathcal{E}}

\newcommand{\sH}{~d\mathcal{H}^{d-1}}

\renewcommand{\emptyset}{\mbox{\textup{\O}}}

\DeclareMathOperator{\divg}{div}

\DeclareMathOperator{\dist}{dist}

\DeclareMathOperator{\Id}{Id}
\DeclareMathOperator{\I}{I}
\DeclareMathOperator{\Err}{Err}

\newcommand{\pD}{\partial D}

\newcommand{\St}{\mathcal{S}}
\newcommand{\Ct}{\mathcal{C}}

\reversemarginpar

\begin{document}

\allowdisplaybreaks

\title[]{A note on the critical set of harmonic functions \\ near the boundary}

\author{Carlos Kenig}
\address{Carlos Kenig
\\ 
Department of Mathematics
\\
University of Chicago
\\
Chicago, IL 60637, USA}
\email{ckenig@uchicago.edu}

\author{Zihui Zhao}
\address{Zihui Zhao
\\ 
Department of Mathematics
\\
Johns Hopkins University
\\
Baltimore, MD 21218, USA}
\email{zhaozh@jhu.edu}

\thanks{The first author was supported in part by NSF grant DMS-2153794, and the second author was partially supported by NSF grant DMS-1902756.}

\date{\today}

\begin{abstract}
	Let $u$ be a harmonic function in a $C^1$ domain $D\subset \RR^d$, which vanishes on an open subset of the boundary. In this note we study its critical set $\mathcal{C}(u): = \{x \in \overline{D}: \nabla u(x) = 0 \}$. When $D$ is a $C^{1,\alpha}$ domain for some $\alpha \in (0,1]$, we give an upper bound on the $(d-2)$-dimensional Hausdorff measure of the critical set by the frequency function. We also discuss possible ways to extend such estimate to all $C^1$-Dini domains, the optimal class of domains for which analogous estimates have been shown to hold for the singular set $\mathcal{S}(u): = \{x \in \overline{D}: u(x) = 0 = |\nabla u(x)| \}$ (see \cite{KZ, KZ-example}).
\end{abstract}

\maketitle

\tableofcontents

\section{Introduction}

The purpose of this note is to study the critical set of harmonic functions that vanish on an open subset of the boundary. In a previous paper \cite{KZ}, we were able to show that for a harmonic function $u$ in a $C^1$-Dini domain $D\subset \RR^d$ and which vanishes on a open set of the boundary, say $B_{2R}(0) \cap \pD$, its singular set $\mathcal{S}(u) := \{x\in \overline{D} \cap B_R(0):  u(x) = 0 = |\nabla u(x)|\}$ is $(d-2)$-rectifiable and its $(d-2)$-dimensional Hausdorff measure is bounded by a constant which depends on the upper bound of the growth rate of $u$ in $B_{2R}(0)$. Here we study the analogous question for the critical set $\mathcal{C}(u) := \{x\in \overline{D} \cap B_R(0): \nabla u(x) = 0 \}$.

We remark that the \L{}ojasiewicz (gradient) inequality (see \cite[(2.4)]{Simon}) implies that for a critical point $x$ of an analytic function $f$, there exists $r>0$ such that all critical points of $f$ in $B_r(x)$ lie on the same level set as $x$, i.e.
\[ \mathcal{C}(f) \cap B_r(x) = \{y\in B_r(x): \nabla f(y) = 0 \} \subset \{f = f(x) \}. \]
In particular, in the analytic setting the critical set and the singular set are essentially the same locally. However, even after a smooth perturbation to the analytic setting, the critical set immediately becomes much more complex than the singular set. For example, let $g: \RR \to \RR$ be an arbitrary smooth function such that $|(g^2)''| = 2|(g')^2 + g g''| < 1/4$. Then $u(x,y,z) : = xy+g^2(z)$ solves the elliptic equation $\divg (A(z) \nabla u) = 0$ in $\RR^3$ with smooth elliptic coefficient matrix
\[ A(z) = \begin{pmatrix}
	1 & -\frac{(g^2)''}{2} & 0 \\
	-\frac{(g^2)''}{2} & 1 & 0 \\
	0 & 0 & 1
\end{pmatrix}. \]
And
\[ \mathcal{C}(u) = \{(0,0,z): g(z) = 0 \text{ or } g'(z) = 0 \}, \quad \mathcal{S}(u) = \{(0,0,z): g(z) = 0\}. \]
Clearly $\mathcal{C}(u)$ can be more complex than $\mathcal{S}(u)$, by choosing suitable smooth functions $g$. (This example was constructed by Leon Simon and appeared in \cite{HHHN}.) Moreover, after a smooth modification of the standard Riemannian metric on the torus $\mathbb{T}^2$, Buhovsky, Logunov and Sodin constructed an infinite sequence of eigenfunctions of the Laplace-Beltrami operator on $(\mathbb{T}^2, g)$, such that each eigenfunction has infinitely many isolated critical points (and the number of singular points is necessarily bounded by the corresponding eigenvalue), see \cite{BLS}.

Nevertheless, in \cite{NVCS} the authors were able to show that for solutions to an elliptic equation with Lipschitz coefficient matrix (without zeroth-order term) in $\RR^d$, one gets $(d-2)$-dimensional volume upper bounds for both $\mathcal{S}(u)$ and $\mathcal{C}(u)$ in the interior (see also \cite{CNV}). The analogous problem at the boundary turns out to be much harder. For examples where (quantitative) unique continuation fails at the boundary, see \cite{BW, Hirsch, Gallegos, KZ-example}. There are serious difficulties in extending our analysis of the singular set in \cite{KZ} to the critical set. In this note, we take a different approach and prove the analogous result for the critical set in $C^{1,\alpha}$ domains, $0<\alpha<1$.
\begin{theorem}\label{thm:main}
	Let $D$ be a $C^{1,\alpha}$ domain in $\RR^d$ with constant $C_\alpha$ ($\alpha\in (0,1)$) such that $0\in \pD$. Let $u$ be a harmonic function in $ D \cap B_{2R}(0)$ such that $u=0$ in $\pD \cap B_{2R}(0)$. Then the critical set
	\[ \mathcal{C}(u) := \{x\in \overline{D}: \nabla u(x) = 0 \} \]
	satisfies that 
	\[ \mathcal{H}^{d-2}(\mathcal{C}(u) \cap B_R(0)) \leq C, \]
	where the constant $C$ depends on the upper bound of the frequency function $N(0,2R)$ of $u$ in $B_{2R}(0)$ (see the definition of the frequency function in \eqref{def:standardfreq}), as well as on $d, R, \alpha$ and $C_\alpha$.
\end{theorem}
\begin{remark}
	\begin{itemize}
		\item In fact, in analogy with the case of the singular set in \cite[Theorem 1.1]{KZ}, the proof gives a stronger estimate on the $(d-2)$-Minkowski content of $\mathcal{C}(u) \cap B_{R}(0)$. Here we state the main theorem for $\mathcal{H}^{d-2}$ for simplicity.
		\item Also, as is the case in \cite[Theorem 1.1]{KZ}, to estimate $\mathcal{C}(u)$ in $B_R(0)$, we need to make assumptions about $u$ in a bigger ball, say $B_{50R}(0)$. Strenuous efforts need to be made to keep track of the enlargement of balls every time, in order to get the stated assumption on $B_{2R}(0)$, but we do not attempt it!
	\end{itemize}
\end{remark}

Our proof of Theorem \ref{thm:main} relies on the following theorem, which is the main result in \cite{HJ}.
\begin{theorem}[Theorem 1.11 in \cite{HJ}]\label{thm:HJ}
	Let $u: B_2 \subset \RR^d \to \RR$ be a non-trivial solution to the elliptic equation $\divg(A\nabla u) = 0$, where the matrix $A(\cdot)$ is H\"older regular with exponent $\alpha \in (0,1)$ and norm $C_\alpha>1$. Assume that $u$ satisfies the following normalized doubling assumption
	\begin{equation}\label{cond:doubling}
		\sup_{B_{2r}(x) \subset B_2} \frac{\mathlarger{\iint}_{B_{2r}(x)} |u-u(x)|^2 ~dy }{\mathlarger{\iint}_{B_{r}(x)} |u-u(x)|^2 ~dy } \leq \Lambda. 
	\end{equation} 
	Then we have 
	\[ \mathcal{H}^{d-2}(\mathcal{C}(u) \cap B_1) \leq C,  \]
	where the constant $C$ depends on $d, \alpha, C_\alpha$ and $\Lambda$.
\end{theorem}
\begin{remark}
\begin{itemize}
	\item We have slightly modified the normalized doubling assumption in \eqref{cond:doubling} to be more consistent with our notation, but it is clearly equivalent to the original assumption in \cite[(1.9)]{HJ}, after modifying the constant $\Lambda$ accordingly. Moreover, since the unique continuation property may fail for solutions to elliptic equations with merely H\"older coefficient matrix (see \cite{Plis, Miller}), this doubling assumption \eqref{cond:doubling} is absolutely necessary to exclude obvious counter-examples.
	\item We again state the conclusion just for the $(d-2)$-dimensional Hausdorff measure for simplicity, but the stronger estimate about the $(d-2)$-dimensional Minkowski content in \cite[Theorem 1.11]{HJ} carries over to Theorem \ref{thm:main}.
\end{itemize}
\end{remark}

We recall that in \cite{KZ}, the analysis for the singular set holds for all $C^1$-Dini domains, a bigger class of domains than $C^{1,\alpha}$ domains (see Definition \ref{def:Dini}). Furthermore, in a follow-up paper \cite{KZ-example}, we show that $C^1$-Dini domain is the necessary assumption. In fact, we construct a large class of \emph{convex} domains which just fail to be $C^1$-Dini domain, such that the $(d-2)$-dimensional Hausdorff measure of the singular set for some (non-negative) harmonic function is infinite.\footnote{See also the work of McCurdy \cite{Mc}, which proves that for \emph{convex domains} the critical set has upper Minkowski dimension $(d-2)$. Our examples show that in general $\mathcal{H}^{d-2}$-finiteness is false, even for the singular set.}
 We remark that our proof of Theorem \ref{thm:main} would hold for all $C^1$-Dini domains, once Theorem \ref{thm:HJ} is proven not just for H\"older-regular coefficient matrices, but for any coefficient matrix that has Dini-regularity, namely
\[ \theta(r): = \sup_{x, y \in B_2 \atop{ |x-y| < r} } |A(y) - A(x)| \text{ satisfies the Dini condition } \int_0^* \frac{\theta(r)}{r} ~dr < \infty. \]
We do not attempt it here.
On the other hand, with a different method the assumption in Theorem \ref{thm:main} can indeed be weakened to $C^1$-Dini domains in the two-dimensional case:

\begin{theorem}\label{thm:twod}
	Let $D$ be a $C^1$-Dini domain in $\RR^2$ with $0\in \pD$. Let $u$ be a harmonic function in $ D \cap B_{2R}(0)$ such that $u=0$ in $\pD \cap B_{2R}(0)$. Then the critical set
	\[ \mathcal{C}(u) := \{x\in \overline{D}: \nabla u(x) = 0 \} \]
	satisfies that 
	\[ \#(\mathcal{C}(u) \cap B_R(0)) \leq C. \]
	where the constant $C$ depends on the Dini parameter of the domain $D$ and $N(0,2R)$, the frequency function of $u$ in the ball $B_{2R}(0)$.
\end{theorem}

In Section \ref{sec:twod}, we prove Theorem \ref{thm:twod} using the tools of conformal mappings.
The proof of Theorem \ref{thm:main} has two ingredients. In Section \ref{sec:domtrsf}, we show that for any $C^1$ domain $D$ and any harmonic function $u$ in $D\cap B_0$ which vanishes on $\pD \cap B_0$, in a smaller ball $B \subset B_0$ we can extend $u$ across the boundary to be a solution to the equation $\divg(A\nabla \tilde{u}) = 0$ in the entire ball $B$, where the modulus of continuity of $A$ is ``one-order less'' than the modulus of continuity of $\pD$. Namely, if locally $D$ is the region above the graph of a $C^1$ function $\varphi$, then the elliptic coefficient matrix $A(\cdot)$ defined on $B$ has a modulus of continuity bounded above by the modulus of continuity of $\nabla \varphi$. In particular, when $D$ is assumed to be a $C^{1,\alpha}$ domain, the corresponding matrix $A(\cdot)$ is $\alpha$-H\"older regular. In Section \ref{sec:doubling}, we study the \emph{normalized frequency function} of $u$ centered at interior points (see Definitions \ref{def:DrHrin} and \ref{def:Nrin}). This has independent interest beyond the goal of this paper, but as an immediate corollary of these results we show that $u$ satisfies the normalized doubling assumption similar to \eqref{cond:doubling}. It follows that the extension of $u$ across the boundary $\tilde{u}$, which is constructed in Section \ref{sec:domtrsf}, also satisfies the normalized doubling assumption. Therefore we can apply Theorem \ref{thm:HJ} to $\tilde{u}$ and estimate its critical set.

Finally, it is somewhat unsatisfactory that Theorem \ref{thm:main} is proven by combining techniques suitable to study \emph{boundary} unique continuations (to verify the normalized doubling assumption \eqref{cond:doubling}) and a reduction to the interior case. It would be more desirable to have a self-contained proof similar to \cite{KZ}, that only uses information at the boundary. In the Appendix, we include some discussions for interested readers about the missing ingredient to make this work. We also give an alternative proof to estimating the critical set in the interior, as a toy model to show why this is not impossible.

\bigskip

\textbf{Acknowledgement.} The second-named author would like to thank Daniele Valtorta for inspiring conversations in the preparation of this work.

\section{Preliminaries}\label{sec:prelim}
\begin{defn}\label{def:Dini}
	Let $\theta: [0,+\infty) \to (0, +\infty)$ be a nondecreasing function verifying
	\begin{equation}\label{cond:Dini}
		\int_0^* \frac{\theta(r)}{r} ~dr < \infty. 
	\end{equation} 
	A connected domain $D$ in $\RR^d$ is a \textit{Dini domain} with parameter $\theta$ if for each point $X_0$ on the boundary of $D$ there is a coordinate system $X=(x,x_d), x\in \RR^{d-1}, x_d\in \RR$ such that with respect to this coordinate system $X_0=(0,0)$, and there are a ball $B$ centered at $X_0$ and a Lipschitz function $\varphi: \RR^{d-1} \to \RR$ verifying the following
	\begin{enumerate}
		\item $\|\nabla \varphi\|_{L^\infty(\RR^{d-1})} \leq C_0$ for some $C_0>0$;
		\item $|\nabla \varphi(x)- \nabla \varphi(y)| \leq \theta(|x-y|)$ for all $x,y \in \RR^{d-1}$;
		\item $D\cap B=\{(x,x_d)\in B: x_d > \varphi(x) \}$.
	\end{enumerate}
	In particular for some exponent $\alpha\in (0,1]$, we say $D$ is a $C^{1,\alpha}$ domain with constant $C_\alpha$, if the parameter $\theta$ satisfies $\theta(r) \leq C_\alpha r^{\alpha}$.
\end{defn}

\begin{defn}
	For any $\theta$ satisfying the assumption in Definition \ref{def:Dini}, we set
\[ \tilde{\theta}(r) = \frac{1}{\log^2 2} \int_r^{2r} \frac1t \int_t^{2t} \frac{\theta(s)}{s} ds ~dt. 
\]
to be a smoothed-out alternative to $\theta$.
Simple computations show that 
\[ \theta(r) \leq \tilde{\theta}(r) \leq \theta(4r), \]
so we can often use $\theta$ and $\tilde{\theta}$ interchangeably, modulo multiplying the variable by a constant.
\end{defn}

For any $p\in \RR^n$, the standard frequency function for $u$ centered at $p$ is defined as
\begin{equation}\label{def:standardfreq}
	N(p,r) := \frac{rD(p,r)}{H(p,r)} = \dfrac{r\iint_{B_r(p) \cap D} |\nabla u|^2 \, dX}{ \int_{\partial B_r(p)\cap D} u^2 \sH}, 
\end{equation} 
c.f. \cite[(6.1) and (6.2)]{KZ}. We remark that it is sometimes more convenient to work with a ``smoothed'' version of the frequency function, namely to define $D(p,r)$ as 
\[ \iint_{B_r(p)\cap D} |\nabla u|^2 \varphi\left( \frac{X-p}{r} \right) \, dX \] 
instead (where $\varphi$ is a smooth cut-off function) and define $H(p,r)$ accordingly (see \cite{DS} and \cite{DMSV}). We use the standard definition for the frequency function to be consistent with \cite{KZ}.
\begin{defn}\label{def:function}
	Let $R, \Lambda>0$ be finite. For any domain $D$ and any function $u$, we say $(u, D) \in \mathfrak{H}(R, \Lambda)$ if
	\begin{itemize}
		\item $D$ is a Dini domain in $\RR^d$ with parameter $\theta$, and it satisfies $\pD \ni 0$, $D$ is graphical inside the ball $B_{5R}(0)$ and
			\begin{equation}\label{cond:R}
				\theta(8R)< \frac{1}{72}, \quad \int_0^{16R} \frac{\theta(s)}{s} ~ds \leq 1 ;
			\end{equation} 
		\item $u$ is a non-trivial harmonic function in $D\cap B_{5R}(0)$,
		\item $u = 0$ on $\pD \cap B_{5R}(0)$,
		\item the \emph{modified} frequency function for $u$ centered at the origin satisfies that $N_0(4R) \leq \Lambda < +\infty$. 
	\end{itemize}  
\end{defn}

In the following we review how we defined the \emph{modified} frequency function centered at boundary points in \cite{KZ}, to guarantee its monotonicity and to have fine control on its derivative. Such fine estimates are not needed for the purpose of this paper, so if the readers do not wish to go into the details, we also add Lemma \ref{lm:starshaped} and some remarks to justify that the last assumption in Definition \ref{def:function} can also be replaced by $N(0,5R) \leq \Lambda < +\infty$, where $N(0,5R)$ is just the standard frequency function as defined in \eqref{def:standardfreq}.


Recall that in \cite[Sections 3 and 4]{KZ} we define the \emph{modified} frequency function centered at boundary points $X\in \pD \cap B_R(0)$ as follows and prove it is monotone increasing with respect to $r$:
\begin{equation}\label{def:modfreq}
	N_X(r) := \widetilde{N}(v_X,r) = N_g(v_X, r) \exp \left(C\int_0^r \frac{\theta(s)}{s} \, ds \right), 
\end{equation} 
where we combine \cite[(4.6) and (3.12)]{KZ}. Here $v_X = u \circ \Psi_X$ with the domain transformation $\Psi_X$ given in \cite[Section 4]{KZ} and $N_g(v_X, r)$ is defined as in \cite[(3.8) and (3.9)]{KZ}. Moreover $N_g(v_X, r)$ is related to the standard frequency function of $u$ as the following lemma indicates.

\begin{lemma}\label{lm:starshaped}
	For any $X\in \pD$ and $r>0$ small (so that $\theta(4r) < 1/26$), we have
	\begin{equation}\label{def:starshaped}
		N_g(v_X, r) = \left[1+O(\theta(4r)) \right] \cdot N(u, X+3r\tilde{\theta}(r) e_d, r), 
	\end{equation} 
	where $N(u, X+3r\tilde{\theta}(r) e_x,r )$ denotes the standard frequency function of $u$ centered at $X+3r \tilde{\theta}(r) e_d \in D$ and with scale $r$, see the definitions in \eqref{def:standardfreq}.
\end{lemma}
\begin{remark}
	This is related to an observation in \cite{KN}: the Dini domain is star-shaped near the boundary. To be more precise, let $X\in \pD$ and $r>0$ be sufficiently small. There exists $Y_r = X+O(r\tilde{\theta}(r))\, n_D(X) \in D$ such that $D \cap B_r(X)$ is star-shaped with respect to $Y_r$. Here $n_D(X)$ denotes the normal derivative of $\pD$ at $X$ pointing inwards. (See the proof of \cite[Lemma 3.2]{KN}.)
\end{remark}
\begin{proof}
	Recall that in \cite[Section 3]{KZ}, we define
	\[ D(v_X, r) = \iint_{B_r \cap \Omega_X} \mu |\nabla_g v_X|^2_g ~dV_g = \iint_{\Psi_{X}(B_r) \cap D} |\nabla u|^2 ~dZ =: \widehat{D}(X, r); \] 
	and by \cite[(3.14) and (4.8)]{KZ}
	\begin{align*}
		H(v_X, r) = \int_{\partial B_r \cap \Omega_X} \mu u^2 ~dV_{\partial B_r}&  = \int_{\partial B_r \cap \Omega_X} \tilde{\eta} v_X^2 \sH \\
		& = \left( 1+ O(\theta(4r)) \right) \int_{\Psi_X(\partial B_r) \cap D} u^2 \sH \\
		& = \left( 1+ O(\theta(4r)) \right) \widehat{H}(X, r),
	\end{align*} 
	where we introduce the notation
	\[ \widehat{H}(X, r) := \int_{\Psi_X(\partial B_r) \cap D} u^2 \sH. 
	\]
	Let
	\[ \widehat{N}(X, r) := \frac{r\widehat{D}(X, r) }{\widehat{H}(X, r) }, \]
	then the frequency function satisfies
	\begin{equation}\label{eq:NvX}
		N(v_X, r) = \frac{r D(v_X, r)}{H(v_X, r)} = \left( 1+O(\theta(4r)) \right) \frac{r\widehat{D}(X, r) }{\widehat{H}(X, r) } = \left( 1+O(\theta(4r)) \right) \widehat{N}(X,r).
	\end{equation}
	
	
	By the definition of the transformation $\Psi$ (see \cite[(4.2)]{KZ}),  we have that
	\[ \Psi(\partial B_r) = \partial B_r + 3r\tilde{\theta}(r) e_d = \partial B_r(3r\tilde{\theta}(r)e_d), \]
	and moreover,
	\[ \Psi(B_r) = B_r(3r\tilde{\theta}(r) e_d), \quad
		\partial \Psi(B_r) = \partial B_r (3r\tilde{\theta}(r) e_d) = \Psi(\partial B_r). \]	
	See the paragraph below \cite[(8.15)]{KZ}.
	Therefore 
	\[ \widehat{H}(X,r) = \int_{\Psi_X(\partial B_r) \cap D} u^2 \sH = \int_{\partial B_r(X+ 3r\tilde{\theta}(r) e_d) \cap D} u^2 \sH, \]
	\[ \widehat{D}(X,r) = \iint_{\Psi_X(B_r) \cap D} |\nabla u|^2 ~dZ = \iint_{B_r(X+ 3r\tilde{\theta}(r) e_d)} |\nabla u|^2 ~dZ, \]
	and the proof is finished.
\end{proof}

Combining \eqref{def:modfreq} and Lemma \ref{lm:starshaped}, we can express the modified frequency function at $X\in \pD$ by the standard frequency function centered at an interior point near $X$:
\[ N_X(r) = \left[1+O(\theta(4r)) \right] \exp \left(C \int_0^r \frac{\theta(s)}{s} \, ds \right) \cdot N(X+3r\tilde{\theta}(r) e_d, r). \]
In particular, an upper bound on the standard frequency function centered at the origin $N(0,5R) \leq \Lambda$ implies an upper bound on the modified frequency function, namely $N_0(4R) \lesssim \Lambda$.

\begin{defn}\label{def:Tpr}
Let $u$ be an arbitrary non-trivial $L^2$ function defined in $D$. We use the following notation to denote its rescalings centered at $X\in \overline{D}$ with scale $r>0$:
\[ T_{X,r}u(Y) := \frac{u(X+r Y)-u(X)}{\left( \frac{1}{r^d} \mathlarger{\iint}_{B_{r}(X)\cap D } \left|u-u(X) \right|^2 ~dZ \right)^{\frac12}}, \qquad T_{r}u(Y) := \frac{u(rY) - u(0)}{\left(\frac{1}{r^d} \mathlarger{\iint}_{B_{r}(0)\cap D } \left|u-u(0)\right|^2 ~dZ \right)^{\frac12}}. \]
%
\end{defn}
\begin{remark}
	Clearly the above rescaling satisfies 
\[ \iint_{B_1(0) \cap \frac{D-X}{r}} |T_{X,r} u|^2 ~dY = 1. \]
\end{remark}

\section{Two-dimensional case}\label{sec:twod}

Since $D$ is graphical, without loss of generality we may assume that $D \cap B_{2R}(0)$ is a simply connected domain bounded by a Jordan curve. If not, we may decompose $D\cap B_{2R}(0)$ into finitely many simply connected domains each containing part of the boundary $\pD$, plus a domain which has some fixed distance away from $\pD$. Then there exists a conformal mapping $\Phi: D\cap B_{2R}(0) \to \RR^2_+$, and moreover by the classical Carath\'eodory theorem (of conformal mapping) $\Phi$ extends to be a homeomorphism from $\overline{D\cap B_{2R}(0)}$ to $\overline{\RR^2_+} \cup \{\infty\}$. Without loss of generality we assume that $\Phi$ maps the origin to the origin, and maps some $X_0 \in \partial B_{2R}(0) \in D$ to $\infty$. Let $g$ denote the imaginary part of $\Phi$, then it is a harmonic function in $D\cap B_{2R}(0)$ which vanishes on $\pD \cap B_{2R}(0)$. Modulo multiplying $\Phi$ by a positive real constant, we may assume that $\sup_{D \cap B_{7R/4}(0)} |g| = 1$. Since $D$ is a $C^1$-Dini domain, by \cite{CK} $g$ is continuously differentiable to the boundary and moreover, its interior normal derivative $\partial_{\nu} g$ has a positive lower bound on $\pD \cap B_{3R/2}(0)$ (Hopf lemma). Because of the Cauchy-Riemann equations, it follows that 
\[ |\det D\Phi | = |\nabla g|^2 = |\partial_\nu g|^2 \geq c>0 \] on $\pD \cap B_{3R/2}(0)$, where $c$ is a universal constant. Moreover, by the continuity of $\nabla g$, we also have $|\det D\Phi| \geq c/2$ when we are sufficiently close to $\partial D \cap B_{3R/2}(0)$ (depending on the modulus of continuity of the Dini parameter of $D$).

Since $\Phi$ maps the origin to the origin, maps $X_0 \in D \cap \partial B_{2R}(0)$ to $\infty$, and it is one-to-one on the boundary, there exists a ball $B'$ containing the origin and centered at a point on $\partial \RR^2_+$, such that $\Phi(\pD \cap B_{3R/2}(0)) = \partial \RR^2_+ \cap B'$. We remark that since both the upper and lower bounds of $|\det D\Phi|$ depend linearly on $\|g\|_{L^{\infty}(B_{7R/4}(0) )}$, which is rescaled to be $1$, the radius of $B'$ satisfies $r_{B'} \approx R$.
 Recall that $u$ is a harmonic function in $D\cap B_{2R}(0)$ which vanishes on $\pD \cap B_{2R}(0)$. Consider $\hat{u}: \RR^2_+ \to \RR$ defined as $\hat{u} := u\circ \Phi^{-1}$.  Since $\Phi$ is a conformal map, $\hat{u}$ is also a harmonic function, and moreover $\hat{u} = 0$ on $\partial \RR^2_+ \cap B'$. Let $B''$ denote a concentric ball as $B'$ such that $r_{B''} \approx r_{B'}$ and $\Phi(D \cap B_R(0)) \supset \RR^2_+ \cap B''$. After an odd reflection and following from a classical result on the interior we have
 \[ \# \left( \mathcal{C}(\hat{u}) \cap B'' \right) \leq C \cdot N_{\hat{u}}(B'),  \]
 where $N_{\hat{u}}(B')$ denotes the standard frequency function for $\hat{u}$ in the ball $B'$.
 
 On the other hand, for any $x \in D\cap B_R(0)$, as long as it is sufficiently close to $\pD \cap B_{3R/2}(0)$ (so that $|\det D\Phi(x)| \geq c/2$), we have that
 \[ x\in \mathcal{C}(u) \iff \Phi(x) \in \mathcal{C}(\hat{u}). \]
 Therefore a simple covering argument gives that
 \[ \#\left(\mathcal{C}(u) \cap B_R(0) \right) \leq C(N_{u}(0,2R)). \]
 

\section{Domain transformation}\label{sec:domtrsf}

Suppose that in some ball $B_0$ centered at the origin, the domain $D$ satisfies
\[ D \cap B_0 = \left\{(x,y) \in \RR^{d-1} \times \RR: y> \varphi(x) \right\} \cap B_0, \]
where $\varphi: \RR^{d-1} \to \RR$ is continuously differentiable. Without loss of generality we assume that $\varphi(0) = 0$ and $\nabla \varphi(0)=0$. (In general when we are studying the local behavior near a point $X_0 \in \pD$ that is not the origin, modulo a translation and rotation we may assume $X_0$ is the origin and that $\pD$ is flat at $X_0$.) Suppose $\theta(\cdot)$ denotes the modulus of continuity of $\nabla \varphi$, namely,
\[ \theta(r) = \sup_{x, z \in \RR^{d-1}, |z|\leq r} |\nabla \varphi(x+z) - \nabla \varphi(x)| \]
In particular, $\varphi$ is a $C^{1,\alpha}$ function (resp. $C^{1,1}$ function, or $C^{1}$-Dini function) if and only if $\theta$ is bounded by $Cr^\alpha$ (resp. $\theta$ is bounded by $Cr$, or $\theta$ satisfies the Dini condition, namely $\int_0^* \frac{\theta(r)}{r} \, dr < \infty$).

Let $X:=(x, \varphi(x))$ be an arbitrary point in $ \pD \cap B_0$. 
Then the unit normal vector to $\pD$ at $X$ pointing inwards is given by
\begin{equation}\label{eq:nv}
	\vec{n}(x) =\left( - \frac{\nabla \varphi(x)}{\sqrt{1+ |\nabla \varphi(x)|^2 } }, \frac{1}{\sqrt{1+|\nabla \varphi(x)|^2}} \right). 
\end{equation} 
In particular, at the origin $\vec{n}(0) = (0,1) = e_d$ points in the upwards vertical direction; and near the origin (which only depends on $\theta$) the modulus of continuity of $\vec{n}$ is controlled by a constant multiple of $\theta$.

Consider the map
\[ G: (x,s) \in \mathbb{R}^{d-1} \times \mathbb{R}_+ \mapsto (x, \varphi(x)) + s \cdot \rho_s \ast \vec{n}(x) \in \mathbb{R}^d, \]
where $\rho \in C_c^\infty(\RR^{d-1})$ is compactly supported in the unit ball centered at the origin, $\rho \geq 0$ and $\int_{\RR^{d-1}} \rho(z) \, dz = 1$, and that $\rho_s(z) := \frac{1}{s^{d-1}} \rho(\frac{z}{s})$ is an approximation to the identity in $\mathbb{R}^{d-1}$. 
In particular, by direct computations, the regularity of $\nabla \varphi$ implies that $G(x,s)$ lies in $D$ when $(x,s)$ is near the origin.
Besides, we have that for any $i=1, \cdots, d-1$,
\begin{align}
	\partial_i G(x,s) & = (e_i, \partial_i \varphi(x)) +  s \cdot \partial_i \left(\rho_s \ast \vec{n}(x) \right) \nonumber \\
	& = (e_i, \partial_i \varphi(x)) +  s \cdot \frac{1}{s} \left( \partial_i \rho\right)_s \ast \vec{n}(x) \nonumber \\
	& = (e_i, \partial_i \varphi(x)) +  \left( \partial_i \rho\right)_s \ast \vec{n}(x);\label{def:DGx}
\end{align} 
and
\begin{align}
	\partial_s G(x,s) & = \rho_s \ast \vec{n}(x) + s \cdot \frac{d}{ds} \left( \rho_s \ast \vec{n}(x) \right) \nonumber \\
	& = \rho_s \ast \vec{n}(x) - s \cdot \frac{1}{s}  \left( (d-1)\rho(z) + z \cdot \nabla \rho(z) \right)_s \ast \vec{n}(x) \nonumber \\
	& = \rho_s \ast \vec{n}(x) - \left( (d-1)\rho(z) + z \cdot \nabla \rho(z) \right)_s \ast \vec{n}(x).\label{def:DGs}
\end{align}

Let $\psi \in C_c^\infty(\RR^{d-1})$ be an arbitrary test function satisfying that $\psi$ is compactly supported in the unit ball and that $\int_{\RR^{d-1}} |\psi(z)| \, dz \leq C$. (For example, we may take $\psi(z) = \rho(z), \partial_i \rho(z)$ or $z\cdot \nabla \rho (z)$.) For any $h\in \RR$ sufficiently small, we have
\begin{align*}
	\left| \psi_{s+h} \ast \vec{n}(x) - \psi_{s} \ast \vec{n}(x) \right| & = \left| \int_{\RR^{d-1}} \vec{n}\left(x-(s+h)z \right) \psi(z) ~dz - \int_{\RR^{d-1}} \vec{n}\left(x-sz \right) \psi(z) ~dz \right| \\
	& \leq \int_{\RR^{d-1}} \left| \vec{n}\left(x-(s+h)z \right) - \vec{n}\left(x-sz \right) \right|  \cdot \left| \psi(z) \right| ~dz \\
	& \lesssim  \int_{\RR^{d-1}} \theta(|hz|)|\psi(z)| ~dz \\
	& \lesssim \theta(|h|),
\end{align*}
where 
in the last inequality we used the assumption that $\psi$ is compactly supported in the unit ball and that $\int_{\RR^{d-1}} |\psi(z)|\, dz$ is bounded.
On the other hand, for any $y\in \RR^{d-1}$ near the origin 
we have
\begin{align*}
	\left| \psi_s \ast \vec{n}(x+y) - \psi_s \ast \vec{n}(x) \right| & = \left|  \int_{\RR^{d-1}} \left[\vec{n}(x+y-z) - \vec{n}(x-z) \right] \psi_s(z) ~dz \right| \lesssim \theta(|y|),
\end{align*}
where the bound is independent of $s$. To summarize, the regularity of $\psi_s \ast \vec{n}(x)$ is controlled by $\theta$. (Away from $s=0$, $\psi_s \ast \vec{n}(x)$ is moreover Lipschitz, but with Lipschitz constant bounded by $C/s$.)
Therefore, by \eqref{def:DGx} and \eqref{def:DGs} the modulus of continuity of $DG(x,s)$, in both the $x$-variable and $s$-variable, is controlled by $\theta$ (uniformly in $s$).


A standard change of variables shows that if $u$ satisfies $\Delta u = 0$ in $D\cap B_0$, then $\tilde{u}(x,s) : = u\circ G(x,s)$ satisfies the equation $-\divg(A(x,s) \nabla \tilde{u}) = 0$ in $\RR^d_+\cap B$, where $A(x,s)$ is a symmetric elliptic matrix given by
\begin{equation}\label{def:Acoefmatrix}
	A(x,s) = |\det DG(x,s)| \left(DG(x,s) \right)^{-1} \left(\left(DG(x,s) \right)^{-1} \right)^T, 
\end{equation} 
and $B$ is a slightly smaller ball than $B_0$, whose radius only depends on $\theta$ (to guarantee that $1/2 \leq \det DG(x,s) \leq 3/2$ in $\RR^d_+ \cap B$).

Notice that 
\[ \int_{\RR^{d-1}} \partial_i \rho(z) \, dz = 0, \]
and
\[ \int_{\RR^{d-1}} z\cdot \nabla \rho(z) \, dz = \sum_{i=1}^{d-1} \int_{\RR^{d-2}} \left( \int_{\RR} z_i \partial_i \rho(z', z_i) \, dz_i \right) d z' = -( d-1) \int_{\RR^{d-1}} \rho(z) \, dz. \]
Therefore, as $s\to 0+$, we have
\[ \partial_i G(x,s) \to (e_i, \partial_i \varphi(x)), \quad \partial_s G(x,s) \to \vec{n}(x), \]
and thus
\[ A(x,s) \to A(x,0) := c(x) \begin{pmatrix}
	\frac{1}{1+|\nabla \varphi(x)|^2 }\left( \Id_{n-1} + \left( \nabla \varphi(x) \right)^T \right) \nabla \varphi(x) & 0 \\
	0 & 1
\end{pmatrix}^{-1}, \]
where 
\[ c(x) = \lim_{s\to 0+} \det DG(x,s) = \frac{1+|\nabla \varphi(x)|^2 }{ \sqrt{1+|\nabla \varphi(x)|^2}} = \sqrt{1+|\nabla \varphi(x)|^2} 
\] is a strictly positive scalar function.

By writing the symmetric matrix $A(x,s)$ on $\RR^d_+\cap B$ as a block matrix
\[ A(x,s) = \begin{pmatrix}
	A_0(x,s) & B(x,s) \\
	B(x,s)^T & d(x,s)
\end{pmatrix}, \]
where $A_0(x,s)$ denotes a $(d-1)\times (d-1)$ matrix,
we define a symmetric matrix $\tilde{A}(x,s)$ on the entire ball $B \subset \RR^d$ as
\[ \tilde{A}(x,s) = \left\{\begin{array}{ll}
	A(x,s), & s > 0 \\
	A(x,0), & s=0 \\
	\begin{pmatrix}
		A_0(x,-s) & -B(x,-s) \\
		-B(x,-s)^T & d(x,-s) \\
	\end{pmatrix}, & s<0.
\end{array} \right. \]
Recall that $\tilde u$ vanishes on $\partial \RR^d_+ \cap B$. Consider its odd reflection across $\RR^d_+$, still denoted as $\tilde{u}$. We have that it satisfies the equation $-\divg(\tilde{A}(x,s) \nabla \tilde{u}) = 0$ in both $\RR^d_+$ and $\RR^d_-$; moreover, since $\lim_{s\to 0\pm} \tilde{A}(x,s) = A(x,0)$ has zeros on the \textit{off-diagonal} block matrix, a simple computation shows that the co-normal derivatives of $\tilde{u}(y,s)$ from above (i.e. $\RR^d_+$) and below (i.e. $\RR^d_-$) cancel each other out, or more precisely,
\[ \lim_{s\to 0+} A(y,s) \nabla \tilde{u}(y,s) \cdot (0,-1) + \lim_{s\to 0-} \tilde{A}(y,s) \nabla \tilde{u}(y,s) \cdot (0,1) = 0. \]
Therefore a standard integration by parts argument shows that $\tilde{u}$ satisfies the equation 
\[ -\divg(\tilde{A}(x,s) \nabla \tilde{u}) = 0 \] in the entire ball $B \subset \RR^d$.
This way, we reduce the study of the boundary value problem for the harmonic function $u$ in $D$ to the study of an interior problem for $\tilde{u}$, which is a solution to an elliptic equation with coefficient matrix $\tilde{A}$. Moreover, by \eqref{def:Acoefmatrix} and the discussion above about the regularity of $DG(x,s)$, the regularity of $\tilde{A}(x,s)$ is also controlled by $\theta$. In particular, if $\Omega$ is a $C^{1,\alpha}$ domain, the matrix $\tilde{A}(x,s)$ is $\alpha$-H\"older regular.

\section{Normalized frequency function and the doubling property}\label{sec:doubling}
In \cite[Section 6]{KZ}, we have studied the standard frequency function for $u$ centered at interior points. In particular, we showed that up to some critical scale the frequency function is monotone (c.f. \cite[Proposition 6.3]{KZ}), and beyond the critical scale the frequency function is always sufficientyl close to the frequency function centered at the corresponding boundary point (c.f. \cite[Lemma 6.17]{KZ}). In this section, we study the \emph{normalized frequency function}, which is better suited to study the critical set.

For any $p\in \RR^d$ and any $r>0$, let
\begin{equation}\label{def:DrHrin}
	D(p, r) = \iint_{B_r(p) \cap D} |\nabla u|^2 ~dX \quad \text{ and } \quad H_{\Ct}(p, r) = \int_{\partial B_r(p) \cap D} |u-u(p)|^2 \sH; 
\end{equation} 
and let 
\begin{equation}\label{def:Nrin}
	N(p,r) = N_{\Ct}(p, r) = \frac{rD(p, r)}{H_{\Ct}(p, r)} 
\end{equation}
whenever the denominator is nontrivial.
Here we use the subindex $\Ct$ since this normalized frequency function is better suited to study the critical set. 
Accordingly, when there is danger of confusion we will add a subindex $\St$ (as in singular set) to the standard frequency function, as is defined in \eqref{def:standardfreq} 

\subsection{Normalized doubling property and the monotonicity formula}

To estimate the normalized frequency function $N(p, r) = N_{\Ct}(p,r)$, we need the following key lemma to compare the $L^2$ norms of $u$ and $u-u(p)$ in the ball $B_r(p)$. Clearly, in general they need not to be comparable in the interior: for example $u$ could have small oscillation but take large value. However, if $B_r(p) \cap \pD \neq \emptyset$, since the value of $u$ is anchored at $0$ at the boundary, this case will not happen and these two norms are indeed comparable.

\begin{lemma}\label{lm:uupsmall}
	Let $R, \Lambda>0$ and $(u, D) \in \HH(R, \Lambda)$. For any $p\in D $ such that $\dist(p, \pD) = \dist(p, \pD \cap B_{\frac{3}{20}R}(0))$, and $r>0$ such that
	\[ \dist(p, \pD) \leq r \leq \frac{3}{16} R, \]
	we have
	\begin{equation}\label{eq:uupsmall}
		\frac{1}{r^{d-1}} \int_{\partial B_r(p) \cap D} |u-u(p)|^2 \sH = \left( 1+ O \left( \left( \frac{\dist(p, \pD)}{r} \right)^{3/4} \right) \right) \frac{1}{r^{d-1}} \int_{\partial B_r(p) \cap D} u^2 \sH.
	\end{equation} 
\end{lemma}
\begin{proof}
	By the sub-harmonicity of $u^2$ and the unique continuation property, it is impossible that
	\[ \int_{\partial B_r(p) \cap D} u^2 \sH = 0, \]
	unless $u \equiv 0$ in $D$. 
	Simple computations show that
	\begin{align*}
		\left| \frac{\frac{1}{r^{d-1}} \mathlarger{\int}_{\partial B_r(p) \cap D} |u-u(p)|^2 \sH }{\frac{1}{r^{d-1}} \mathlarger{\int}_{\partial B_r(p) \cap D} u^2 \sH } - 1 \right| & = \left|\frac{-2u(p) \cdot \frac{1}{r^{d-1}} \mathlarger{\int}_{\partial B_r(p) \cap D} u \sH + u^2(p) \frac{\mathcal{H}^{d-1}(\partial B_r(p) \cap D)}{r^{d-1}} }{\frac{1}{r^{d-1}} \mathlarger{\int}_{\partial B_r(p) \cap D} u^2 \sH } \right| \\
		& \lesssim \frac{|u(p)|}{ \left(\frac{1}{r^{d-1}} \mathlarger{\int}_{\partial B_r(p) \cap D} u^2 \sH \right)^{1/2}} + \frac{|u(p)|^2}{\frac{1}{r^{d-1}} \mathlarger{\int}_{\partial B_r(p) \cap D} u^2 \sH}.
	\end{align*}
	Hence it suffices to estimate the ratio and show
	\begin{equation}\label{eq:ptest}
		\frac{|u(p)|}{ \left(\frac{1}{r^{d-1}} \mathlarger{\int}_{\partial B_r(p) \cap D} u^2 \sH \right)^{1/2}} \lesssim \left( \frac{\dist(p,\pD)}{r} \right)^{3/4}. 
	\end{equation} 
	
	Let $q\in \pD \cap B_{\frac{3}{20}R}(0)$ such that $|p-q| = \dist(p,\pD)$. Set $s=\dist(p,\pD)$ for simplicity of notation.
	Since $u$ vanishes on the boundary, \eqref{eq:ptest} follows from the boundary H\"older regularity of $u$ with H\"older exponent sufficiently close to $1$ in $C^1$ domains (and the sub-harmonicity of $u^2$), see for example \cite{Dahlberg}.
	We give an alternative proof here in order to be self-contained. 
	 By the boundary Harnack inequality (or the sub-harmonicity of $u^2$), we have
	\begin{align*}
		|u(p)|^2 \lesssim \frac{1}{s^{d}} \iint_{B_{2s}(q) \cap D} u^2 ~dZ. 
	\end{align*}
	On the other hand, since 
	$2s \leq 2r$, by the doubling property (see \cite[(3.29)]{KZ}) we have 
	\begin{align*}
		\frac{\frac{1}{(2r)^{d}} \mathlarger{\iint}_{B_{2r}(q) \cap D} u^2 \sH }{\frac{1}{(2s)^{d}} \mathlarger{\int}_{B_{2s}(q) \cap D} u^2 \sH} \gtrsim \left( \frac{r}{s} \right)^{2N_q \cdot\exp\left(-C\int_0^{2r} \frac{\theta(\tau)}{\tau} \, d\tau \right)} \geq \left( \frac{r}{s} \right)^{3/2},
	\end{align*}
	where we also used the monotonicity of $N_q(\cdot)$, the fact that $N_q = \lim_{\rho \to 0} N_q(\rho) \geq 1$ (see \cite[Lemma 5.23]{KZ}), and the smallness of $\int_0^{2r} \frac{\theta(\tau)}{\tau} \, d\tau$ for $r\leq \frac{3}{16} R$ sufficiently small.
	Therefore
	\begin{align}
		|u(p)|^2 \lesssim \left( \frac{s}{r} \right)^{3/2} \cdot  \frac{1}{(2r)^{d}} \iint_{ B_{2r}(q) \cap D} u^2 ~dZ & \lesssim \left( \frac{s}{r} \right)^{3/2} \cdot  \frac{1}{(3r)^{d}} \iint_{B_{3r}(p) \cap D} u^2 ~dZ \nonumber \\
		& \lesssim_\Lambda \left( \frac{s}{r} \right)^{3/2} \cdot  \frac{1}{r^{d}} \iint_{B_{r}(p) \cap D} u^2 ~dZ.\label{tmp:uestdist}
	\end{align}
	In the last inequality we use $N_{\St}(p,3r) \leq C(\Lambda)$, which follows easily by comparing $N_{\St}(p,3r)$ and $N_q(4r)$ using Lemma \ref{lm:starshaped} and the upper bound of $N_q(4r)$ (c.f. \cite[(5.2)]{KZ}). Again \eqref{tmp:uestdist} implies \eqref{eq:ptest} by the sub-harmonicity of $u^2$.
\end{proof}

\begin{lemma}\label{lm:uup}
	Let $R, \Lambda>0$ be fixed and $(u,D) \in \mathfrak{H}(R, \Lambda)$. There exist constants $c_2\geq c_1>0$ such that for any $p\in D\cap B_{\frac{3R}{40}}(0)$ and $\dist(p, \pD) \leq r< \frac{R}{8}$, we have
	\begin{equation}\label{eq:uSuC}
		c_1 \frac{1}{r^d} \iint_{B_r(\tilde{p})} u^2 ~dZ \leq \frac{1}{r^{d}} \iint_{B_r(p) \cap D}|u-u(p)|^2 ~dZ \leq c_2 \frac{1}{r^d} \iint_{B_r(\tilde{p})} u^2 ~dZ, 
	\end{equation} 
	where $\tilde{p} \in \pD$ satisfies $|\tilde{p} - p| = \dist(p, \pD)$.
\end{lemma}
\begin{proof}
	The upper bound follows easily from Lemma \ref{lm:uupsmall}, but the lower bound follows only when $\dist(p, \pD) \ll r$. So instead we prove the lower bound by contradiction.
	Assume not, then there exist $(u_i, D_i) \in \mathfrak{H}(R, \Lambda)$ and $p_i \in D_i \cap B_{\frac{3R}{40}}(0)$ satisfying $\dist(p_i, \pD_i) \leq r_i < \frac{R}{8}$, such that
	\begin{equation}\label{cl:tmpuup}
		\frac{1}{r_i^{d}} \iint_{B_{r_i}(p_i) \cap D_i}|u_i-u_i(p_i)|^2 ~dZ < 2^{-i} \frac{1}{r_i^d} \iint_{B_{r_i}(\tilde{p}_i)} u_i^2 ~dZ 
	\end{equation} 
	where $\tilde{p}_i\in \pD_i$ satisfies $|\tilde{p}_i -p_i| = \dist(p_i, \pD_i)$.
	Since $\iint_{B_{r_i}(\tilde{p}_i)} u_i^2 ~dZ \neq 0$ when $u_i$ is assumed to be a nontrivial harmonic function, we may assume without loss of generality that
	\[ \frac{1}{r_i^d} \iint_{B_{r_i}(\tilde{p}_i)} u_i^2 ~dZ = 1. \]
	Thus \eqref{cl:tmpuup} is simplified to 
	\begin{equation}\label{cl:uup}
		\frac{1}{r_i^{d}} \iint_{\partial B_{r_i}(p_i) \cap D_i}|u_i-u_i(p_i)|^2 ~dZ < 2^{-i}. 
	\end{equation}
	
	By the compactness of rescaled harmonic functions at boundary points (see \cite[Proposition 5.24 and Remark 5.25]{KZ}), the sequence $T_{\tilde{p}_i, r_i} u_i$ converges to a harmonic function $u_\infty$ in $D_\infty \cap B_5(0)$, where
	\[ D_\infty = \lim_{i\to \infty} \frac{D_i-\tilde{p}_i}{r_i} \text{ is a connected Lipschitz domain with small Lipschitz constant in } B_5(0) \]
	and $\pD_\infty \ni 0$.
	Besides, since $|p_i- \tilde{p}_i| = \dist(p_i, \pD_i) \leq r_i$, modulo passing to a subsequence $q_i:=\frac{p_i-\tilde{p}_i}{r_i}$ converges to some point $q\in \overline{D_\infty} \cap \overline{B_1(0)}$. We can rewrite \eqref{cl:uup} using the rescaled functions $T_{\tilde{p}_i, r_i} u_i(Z) = u_i(\tilde{p}_i + r_i Z)$ as
	\[ \iint_{B_1(q_i) \cap \frac{D_i- \tilde{p}_i}{r_i}} \left|T_{\tilde{p}_i, r_i} u_i(Z) - T_{\tilde{p}_i, r_i}u_i(q_i) \right|^2 ~dZ < 2^{-i}. \]
	By the uniform convergence of $T_{\tilde{p}_i, r_i}u_i$ to $u_\infty$, graphical convergence of the domain $\frac{D_i- \tilde{p}_i}{r_i}$ to $D_\infty$, $q_i \to q$ and the uniform gradient bound of $T_{\tilde{p}_i, r_i}u_i$ (by \cite[Lemma 5.19]{KZ}), we may pass to the limit and conclude
	\begin{equation}\label{eq:uinftyconst}
		\iint_{B_{1/2}(q) \cap D_\infty} |u_\infty - u_\infty(q)|^2 ~dZ = 0. 
	\end{equation} 
	Hence 
	\begin{equation}\label{tmp:constonsph}
		u_\infty \equiv u_\infty(q) \text{ on } B_{1/2}(q) \cap D_\infty.
	\end{equation}
	Combined with the unique continuation property in the connected Lipschitz domain $D_\infty \cap B_5(0)$ and the fact that $u_\infty$ has vanishing boundary data, we get a contradiction with the fact that $\iint_{B_1(0)} u_\infty^2 ~dZ = 1$.
	
\end{proof}
\begin{remark}\label{rm:lmuup}
	By an analogous proof, one can also show that \eqref{eq:uSuC} holds when the middle term is replaced by 
	\[ \frac{1}{r^{d-1}} \int_{\partial B_r(p) \cap D} |u-u(p)|^2 \sH. \]
	We explain the necessary modification of the proof below. By an analogous argument by contradiction, we obtain in lieu of \eqref{eq:uinftyconst} the following
	\[ \int_{\partial B_1(q) \cap D_\infty} |u_\infty - u_\infty(q)|^2 \sH = 0. \]
	Hence 
	\begin{equation}\label{tmp:constonsphere}
		u_\infty \equiv u_\infty(q) \text{ on } \partial B_1(q) \cap D_\infty. 
	\end{equation} 
	Recall that $D_\infty$ is a Lipschitz domain with small Lipschitz constant in $B_5(0)$, $0\in \partial D_\infty$ and $q\in \overline{D_\infty} \cap \overline{B_1(0)}$. In particular $\partial B_1(q) \cap \partial D_\infty \neq \emptyset$. This, combined with \eqref{tmp:constonsphere} and the fact that $u_\infty$ vanishes on $\partial D_\infty \cap B_5(0)$, implies that $u_\infty \equiv 0$ on $\partial (B_1(q) \cap D_\infty)$. By the maximum principle $u_\infty \equiv 0$ on $B_1(q) \cap D_\infty$, and thus $u_\infty \equiv 0$ on all of $D_\infty$ (at least inside $B_5(0)$ when it is explicitly known to be connected) by the unique continuation property. This again contradicts with the fact that $\iint_{B_1(0)} u_\infty^2 ~dZ = 1$.
\end{remark}

\begin{prop}\label{prop:intmonotonicity}
	Let $D$ be a Dini domain, 
	and $u$ be a harmonic function in $D$ such that $u =0 $ on $\pD \cap B_{5R}(0)$. Then for any $p\in D\cap B_{\frac{3R}{40}}(0)$, the normalized frequency function $N(p, r) = N_{\Ct}(p,r)$ satisfies
	\begin{equation}\label{eq:dernNr}
		N'(r) = \left\{\begin{array}{ll}
			R_h(r) + R_b(r), & \text{ if } r< \dist(p,\pD) \\
			R_h(r) + R_b(r) +  O\left( \left( \frac{\dist(p, \pD)}{r} \right)^{3/4} \right) \cdot \frac{N(r)}{r} , & \text{ if } r \geq \dist(p, \pD)
		\end{array} \right.  
	\end{equation} 
	where
	\begin{equation}\label{def:Rhr}
		R_h(r):= \frac{2r}{H(r)} \int_{\partial B_r(p) \cap D} \left|\partial_\rho u - \frac{N(r)}{r} \left(u-u(p) \right)  \right|^2 \sH, 
	\end{equation} 
	\[ R_b(r):= \frac{1}{H(r)} \int_{B_r(p) \cap \pD} (\partial_n u)^2 \langle X-p, n_D(X) \rangle \sH. \]
	Here $\partial_\rho u$ denotes the radial derivative of $u$ with vertex $p$, and $\partial_n u$ denotes the normal derivative of $u$ pointing away from $D$.
	
	Moreover $R_b(r) \geq 0$
	as long as 
	\begin{equation}\label{cond:intmono}
		r\theta(r) \leq \dist(p,\pD).
	\end{equation} 
\end{prop}
\begin{proof}
	Since the proof is similar to that of \cite[Proposition 6.3]{KZ}, we just highlight the differences here. Since the definition of $D(r)$ has not changed, the derivative of $D(r)$ stays the same:
	\begin{equation}\label{eq:nderDr}
		D'(r) = \frac{d-2}{r} D(r) + 2\int_{\partial B_r(p) \cap D} (\partial_\rho u)^2 \sH + \frac1r \int_{B_r(p) \cap \pD} (\partial_n u)^2 \langle X-p, n_D(X) \rangle \sH.
	\end{equation}
	However, since $H(r) = H_{\Ct}(r)$ is now defined as the $L^2$ norm of $u-u(p)$, we have
	\begin{equation}\label{eq:nderHr}
		H'(r) = \frac{d-1}{r} H(r) + 2 \int_{\partial B_r(p) \cap D} \left(u-u(p) \right) ~\partial_\rho u \sH.
	\end{equation}
	Again by the divergence theorem and the vanishing boundary condition,
	\begin{align}
		D(r) & = \iint_{B_r(p) \cap D} |\nabla(u-u(p))|^2 ~dX \nonumber \\
		& = \iint_{B_r(p) \cap D} \divg\left((u-u(p)) \nabla u \right) ~dX \nonumber \\
		& = \int_{\partial B_r(p) \cap D} \left( u-u(p) \right) ~\partial_{\rho} u\sH + \int_{B_r(p) \cap \pD} \left(u - u(p) \right) ~\partial_n u \sH \nonumber \\
		& = \int_{\partial B_r(p) \cap D} \left( u-u(p) \right) ~\partial_{\rho} u\sH - u(p) \cdot \int_{B_r(p) \cap \pD} \partial_n u \sH.\label{eq:nderHrDr}
	\end{align}
	Therefore combining \eqref{eq:nderDr}, \eqref{eq:nderHr} and \eqref{eq:nderHrDr}, we get
	\begin{align*}
		\frac{N'(r)}{N(r)} & = \frac1r + \frac{D'(r)}{D(r)} - \frac{H'(r)}{H(r)} \\
		& = \frac{2}{D(r)} \int_{\partial B_r(p) \cap D} \left|\partial_\rho u - \frac{N(r)}{r} \left(u-u(p) \right)  \right|^2 \sH \\
		& \qquad + \frac{1}{rD(r)} \int_{B_r(p) \cap \pD} (\partial_n u)^2 \langle X-p, n_D(X) \rangle \sH \\
		& \qquad \qquad + \frac{2u(p)\cdot \int_{B_r(p) \cap \pD} \partial_n u\sH }{H(r)}.
	\end{align*}
	It suffices to estimate the error term
	\begin{equation}\label{def:Err}
		\Err_r: = \frac{2u(p)\cdot \int_{B_r(p) \cap \pD} \partial_n u\sH }{H(r)}. 
	\end{equation} 
	
	We first remark that this error term only appears if $r\geq \dist(p,\pD)$, so the first part of \eqref{eq:dernNr} is proven. Let $\tilde{p} \in \pD$ satisfy $|\tilde{p}-p| = \dist(p, \pD)$. Recall we have the following gradient estimate for harmonic functions with vanishing boundary data in Dini domains (see for example \cite{CK}):
	\[ \sup_{B_{r}(\tilde{p}) \cap D} r|\nabla u| \lesssim \left(\frac{1}{r^d} \iint_{B_{2r}(\tilde{p}) \cap D} u^2 ~dX \right)^{1/2}. \]
	On the other hand, we have shown in \eqref{tmp:uestdist} that
	\[ |u(p)| \lesssim \left( \frac{\dist(p,\pD)}{r} \right)^{3/4} \left(\frac{1}{r^d} \iint_{B_{2r}(\tilde{p}) \cap D} u^2 ~dX \right)^{1/2}. \]
	Therefore
	\begin{align*}
		|\Err_r| \leq \frac{2}{H(r)} |u(p)|\cdot \int_{B_r(p) \cap \pD} |\nabla u| \sH & \lesssim \frac{1}{H(r)}\cdot \left( \frac{\dist(p,\pD)}{r} \right)^{3/4} \cdot \frac{1}{r^2} \iint_{B_{2r}(\tilde{p}) \cap D} u^2 ~dX \\
		& \lesssim_{\Lambda} \left( \frac{\dist(p,\pD)}{r} \right)^{3/4} \cdot \frac{1}{r^2H(r)}  \iint_{B_r(\tilde{p}) \cap D} u^2 ~dX \\
		& \lesssim \left( \frac{\dist(p,\pD)}{r} \right)^{3/4} \frac{1}{r} \\
		& =: C_m \left( \frac{\dist(p,\pD)}{r} \right)^{3/4} \frac{1}{r},
	\end{align*}
	where we use the $L^2$-doubling property centered at $\tilde{p}$ in the penultimate inequality, and we use Lemma \ref{lm:uup} and Remark \ref{rm:lmuup} in the last inequality. We denote $C_m$ as the multiplication of all constants appearing in the inequalities. This proves \eqref{eq:dernNr}.
	
	The last part follows because the condition \eqref{cond:intmono} implies that
	\[ \langle X-p, n_D(X) \rangle \geq 0 \quad \text{ for every } X\in B_r(p) \cap \pD. \]
	This is proved in \cite[Proposition 6.3]{KZ}. 
\end{proof}

As in \cite[Section 6.1]{KZ}, for any any $p\in D\cap B_{2R}(0)$ sufficiently close to the boundary, we refer to the unique $r>0$ such that \begin{equation}\label{def:criticalval}
	 \dist(p, \pD) = r\tilde{\theta}(r) 
\end{equation} as the \emph{critical scale} for $p$, and denote it as $r_{cs}(p)$. In particular, by the assumption on $\theta(R)$ in \eqref{cond:R}, we always have $r_{cs}(p) \gg \dist(p,\pD)$ (to be more precise $r_{cs}(p) \geq 72 \dist(p,\pD)$). 

It follows easily from Proposition \eqref{prop:intmonotonicity} that $N(r)$ is monotone non-decreasing when $0<r<\dist(p,\pD)$. For any radii $r_1 < r_2$ within the interval $[\dist(p,\pD), r_{cs}(p)]$, since
\[ \frac{N'(r)}{N(r)} \geq -C_m \left( \frac{\dist(p,\pD)}{r} \right)^{3/4} \frac{1}{r}, \]
we have that
\begin{align*}
	\log\frac{N(r_2)}{N(r_1)} = \int_{r_1}^{r_2} \frac{N'(r)}{N(r)} ~dr  & \geq -C_m \int_{r_1}^{r_2} \left( \frac{\dist(p,\pD)}{r} \right)^{3/4} \frac{dr}{r} \\
	& \geq -\frac43 C_m \left( \frac{\dist(p,\pD)}{r_1} \right)^{3/4} \geq - \frac43 C_m.
\end{align*}
Hence
\begin{equation}\label{eq:monmodconst}
	N(r_1) \leq \exp\left( \frac43 C_m \right) N(r_2), \quad \text{ for any } \dist(p,\pD) \leq r_1 < r_2 \leq r_{cs}(p).
\end{equation} 
In particular, though we may lost monotonicity on $[\dist(p,\pD), r_{cs}(p)]$, we do still have the \emph{boundedness} of the normalized frequency function.

As in \cite[Lemma 5.1]{KZ}, we also have a rough uniform bound (depending on $d, R, \Lambda$) on the normalized frequency function.
\begin{lemma}\label{lm:freqbd_in}
	Assume that $N_0(4R) \leq \Lambda<+\infty$. Suppose $X\in D \cap B_{\frac{R}{40}}(0)$ is such that $\dist(X, \pD) \leq \frac{R}{8}$. Then for any $r\leq  \frac{R}{8}$, we have
	\begin{equation}\label{eq:freqbdR_in}
		N_{\mathcal{C}}(X, r) \leq C(\Lambda) N_{\tilde{X}} \left( \frac34 R \right),
	\end{equation} 
	where $\tilde{X} \in \pD$ satisfies $|X-\tilde{X}| = \dist(X, \pD)$.
%
%
\end{lemma}
\begin{remark}
	Notice that the above upper bound holds even outside of the monotonic interval $[0,r_{cs}(X)]$, as long as we have $r\leq \frac{R}{8}$.
\end{remark}
\begin{proof}
	Since the normalized frequency function $N(X, \cdot) = N_{\mathcal{C}}(X,\cdot)$ is monotone non-increasing (up to a constant, see \eqref{eq:monmodconst}) on the interval $[0, r_{cs}(X)]$ and $\dist(X, \pD) < r_{cs}(X)$, it suffices to prove \eqref{eq:freqbdR_in} for $\dist(X, \pD) \leq r \leq \frac{R}{8}$.
	Recall that we have proven in Lemma \ref{lm:uup} and Remark \ref{rm:lmuup} that
	\[ H(X,r) \gtrsim \frac{1}{r} \iint_{B_r(\tilde{X})} u^2 ~dZ. \]
	On the other hand by the uniform bound on $N_{\tilde{X}}(r)$ and the doubling property in \cite[(3.30)]{KZ} we have
	\begin{align*}
		\frac{1}{r} \iint_{B_r(\tilde{X})} u^2 ~dZ = \frac{1}{r}  \int_0^r H(\tilde{X},s) ~ds \gtrsim \frac{1}{r} \int_0^r H(\tilde{X},2r) \left( \frac{s}{2r} \right)^{d-1+C(\Lambda)} ~ds \gtrsim_{d,\Lambda} H(\tilde{X}, 2r).
	\end{align*}
	Therefore
	\begin{align*}
		N(X, r) = \frac{r \mathlarger{\iint}_{B_r(X)} |\nabla u|^2 ~dZ}{ H(X, r)} \lesssim \frac{r \mathlarger{\iint}_{B_{2r}(\tilde{X}) } |\nabla u|^2 ~dZ }{H(\tilde{X}, 2r)} \leq N_{\tilde{X}}(2r) \leq N_{\tilde{X}}\left(\frac34 R \right).
	\end{align*}
%
%
%
\end{proof}

The above frequency upper bound yields the doubling property for the $L^2$ norm of $u-u(X)$. 
\begin{corollary}\label{cor:doublingu_in}
	Assume that $N_0(4R) \leq \Lambda<+\infty$. Then for any $X\in D\cap B_{\frac{R}{40}}(0)$ such that $\dist(X, \pD) \leq \frac{R}{8}$,
	and any $\rho>0$, $1<a<72$ such that $a\rho \leq \frac{R}{8}$, we have
	\begin{equation}\label{eq:doublingin}
		\frac{\mathlarger{\iint}_{B_{a\rho}(X)\cap D} |u-u(X)|^2 ~dZ}{\mathlarger{\iint}_{B_{\rho}(X) \cap D} |u-u(X)|^2 ~dZ} \lesssim a^{d+C(\Lambda,R)}.
	\end{equation}
\end{corollary}
\begin{proof}
Recall that on the interval $[0, r_{cs}(X)]$, we have computed the derivative of $H_{\Ct}(X,s)$ in \eqref{eq:nderHr} and \eqref{eq:nderHrDr}, which can be reformulated into
\begin{align*}
	\left(\log \frac{H_{\Ct}(X,s)}{s^{d-1}} \right)' = \frac{2N_{\Ct}(X,s)}{s} + \Err_s,
\end{align*}
where $\Err_s$ is defined as in \eqref{def:Err} and it satisfies $|\Err_s| \leq C_m/s$. Integrating the above equality and by the upper bound of the frequency function in Lemma \ref{lm:freqbd_in}, we get
\begin{equation}\label{eq:doublingnH}
	\frac{H_{\Ct}(X,as)}{H_{\Ct}(X,s)} = a^{d-1}\exp\left[ \int_s^{as} \left(\frac{2N_{\Ct}(X,\tau)}{\tau} +  \Err_{\tau} \right) ~d\tau \right] \leq a^{d-1+C(\Lambda, R)}. 
\end{equation} 
Therefore
\begin{align}
	\frac{\mathlarger{\iint}_{B_{a\rho}(X) \cap D} |u-u(X)|^2 ~dZ}{\mathlarger{\iint}_{B_{\rho}(X) \cap D} |u-u(X)|^2 ~dZ} = \frac{a\mathlarger{\int}_0^{\rho} H_{\Ct}(X,as) ~ds}{\mathlarger{\int}_0^{\rho} H_{\Ct}(X,s) ~ds } \leq a^{d+C(\Lambda,R)} \label{tmp:doubling1}
\end{align}
as long as $a\rho \leq r_{cs}(X)$.

On the other hand, when  $\dist(X,\pD) \leq \rho < a\rho < \frac{R}{8}$, by Lemma \ref{lm:uup} we have
\begin{equation}\label{tmp:doubling2}
	\frac{\mathlarger{\iint}_{B_{a\rho}(X) \cap D}|u-u(X)|^2 ~dZ }{\mathlarger{\iint}_{B_{\rho}(X) \cap D}|u-u(X)|^2 ~dZ} \lesssim \frac{\mathlarger{\iint}_{B_{a\rho}(\tilde{X})} u^2 ~dZ }{\mathlarger{\iint}_{B_{\rho}(\tilde{X})} u^2 ~dZ} \lesssim a^{d+C(\Lambda, R)}, 
\end{equation} 
where we used the upper bound of the frequency function $N_{\tilde{X}}(\cdot)$ and the doubling property at boundary points, see for example \cite[Corollary 3.28]{KZ}.

Since $r_{cs}(X) \geq 72 \dist(X, \pD)$, we can combine the estimates \eqref{tmp:doubling1} and \eqref{tmp:doubling2} to get the final conclusion whenever $a\rho < \frac{R}{8}$.
\end{proof}

\subsection{Normalized frequency function beyond the critical scale}

The normalized frequency function $N_{\mathcal{C}}(X, \cdot)$ is only controlled on the interval $[0,r_{cs}(X)]$; for large radius, it ought to be replaced by the corresponding frequency function centered at a boundary point $\tilde{X} \in \pD$, which satisfies $|X-\tilde{X}| = \dist(X, \pD)$, as in the case for the standard frequency function $N_{\mathcal{S}}(X, \cdot)$ (c.f. \cite[Lemma 6.17]{KZ}).

First we need to the following estimate, which is analogous to \cite[Lemma 6.14]{KZ}.

\begin{lemma}\label{lm:urunifgrad}
	Assume that $N_0(4R)\leq \Lambda +\infty$. For any $X \in D \cap B_{\frac{R}{40}}(0)$ such that $\dist(X,\pD)\leq \frac{R}{8}$, and any
	$r\leq \frac{R}{48}$, we have
	\[ \left|\nabla T_{X, r} u(Y) \right| \leq C, \quad \text{ for all } Y \in \frac{D-X}{r} \cap B_5. \]
\end{lemma}
\begin{proof}
	Recalling the definition of $T_{X,r}u$ in Definition \ref{def:Tpr}, we have
	\[ |\nabla T_{X,r}u(Y)|^2 = \frac{\left( r|\nabla u(X+rY)| \right)^2}{\frac{1}{r^d} \mathlarger{\iint}_{B_r(X)\cap D} |u-u(X)|^2 ~dZ }. \]
	We first consider the purely interior case when $6r < \dist(X, \pD)$.
	In this case, by the interior gradient estimate and the doubling property in Corollary \ref{cor:doublingu_in}, we have
	\begin{align*}
		\left|\nabla T_{X, r} u(Y) \right| \leq \frac{\sup_{B_{5r}(X)} \left( r|\nabla u| \right)^2}{\frac{1}{r^d} \mathlarger{\iint}_{B_r(X)\cap D} |u-u(X)|^2 ~dZ } \lesssim \frac{ \mathlarger{\iint}_{B_{6r}(X)\cap D} |u-u(X)|^2 ~dZ}{\mathlarger{\iint}_{B_r(X)\cap D} |u-u(X)|^2 ~dZ } \leq 6^{d+C(\Lambda,R)}
	\end{align*}
	is uniformly bounded. Now assume $\dist(X, \pD) \leq 6r \leq \frac{R}{8}$.	
	Let $\tilde{X}\in \pD$ such that $|X-\tilde{X}| = \dist(X,\pD)$.
	We bound the numerator by the gradient estimate of $u$:
	\begin{equation}\label{tmp:numbd}
		\left( r|\nabla u(X+rY)| \right)^2 \leq \sup_{B_{6r}(\tilde{X}) \cap D} \left( r|\nabla u|\right)^2 \lesssim \frac{1}{r^d} \iint_{B_{7r}(\tilde{X})} u^2 ~dZ.
	\end{equation}
	To bound the denominator, we use the doubling property \eqref{eq:doublingin} and Lemma \ref{lm:uup} to get
	\begin{align}
		\frac{1}{r^d} \iint_{B_r(X)\cap D} |u-u(X)|^2 ~dZ \gtrsim \frac{1}{(6r)^d} \iint_{B_{6r}(X) \cap D} |u-u(X)|^2 ~dZ \gtrsim_{c_1} \frac{1}{(6r)^d} \iint_{B_{6r}(\tilde{X})} u^2 ~dZ. \label{tmp:denbd}
	\end{align}
	Therefore combining \eqref{tmp:numbd} and \eqref{tmp:denbd}, it follows from \cite[(3.29)]{KZ} 
	that
	\[ |\nabla T_{X,r} u|^2 \lesssim \frac{\mathlarger{\iint}_{B_{7r}(\tilde{X})} u^2 ~dZ }{\mathlarger{\iint}_{B_{6r}(\tilde{X})} u^2 ~dZ } \lesssim \left(\frac76 \right)^{d+N_{\tilde{X}}(7r)} \leq C(\Lambda). \]
	
%
%
\end{proof}


\begin{lemma}\label{lm:spvarin_far}
	Let $R, \Lambda>0$, $\rho \in (0, 1/6]$ and $\delta_{in}>0$ be fixed. There exists $r_{in} = r_{in}(\delta_{in}, \rho)>0$ such that the following holds for any $(u,D) \in \HH(R,\Lambda)$. Suppose $p\in D\cap B_{\frac{R}{40}}(0)$.
	Let $q\in \pD $ satisfy $|q-p| \leq 2 \dist(p,\pD)$. Then for any radius $r$ such that $\rho r_{cs}(p) \leq r \leq r_{in}$
	\footnote{Note that for such radius, the normalized frequency function $N_{\mathcal{C}}(p, r)$ is no longer monotone; but it is still defined.} 
	we have
	\[ \left| N_{\mathcal{C}}(p, r) - N_q( r) \right| \leq \delta_{in}. \]
%
\end{lemma}

\begin{proof}
	We argue by contradiction. To that end we assume there exist sequences $(u_i, D_i) \in \HH(R, \Lambda)$, $p_i \in D_i\cap B_{\frac{R}{40}}(0)$ with $\rho r_{cs}(p_i) \leq r_i \to 0$ and 
	\begin{equation}\label{assp:distr}
		d_{p_i} := \dist(p_i, \pD_i) = r_{cs}(p_i) \tilde{\theta}(r_{cs}(p_i)) \leq \frac{r_i}{\rho} \tilde{\theta}\left(\frac{r_i}{\rho} \right), 
	\end{equation}  
	and $q_i \in \pD_i$ with $|p_i - q_i| \leq 2d_{p_i}$, such that
	\begin{equation}\label{cl:freqbin}
		\left|N_{\mathcal{C}}(p_i, r_i) - N_{q_i}(r_i) \right| > \delta_{in}>0. 
	\end{equation}
	Recall that $N_{q_i}(r_i) = \widetilde{N}(v_{q_i}, r_i)$, see \eqref{def:modfreq} as well as the definition that $v_{q_i} = u \circ \Psi_{q_i}$ in \cite[Section 3]{KZ}.
	
	By assumption 
	\[ D_i \cap B_{5R}(0) \subset \{(x, x_d) \in \RR^{d-1} \times \RR: x_d > \varphi_i(x)\} \]
	for some $C^1$ function $\varphi_i$ with Dini parameter $\theta$.
	Without loss of generality we assume $\varphi_i(0) = 0$ and $\nabla \varphi_i(0) = 0$. Hence 
	\[ |\nabla \varphi_i(x)| = |\nabla \varphi_i(x) - \nabla \varphi_i(0) | \leq \theta(|x|) \leq \theta(5R) \]
	is uniformly bounded, and by Arzela-Ascoli $\varphi_i$ converges uniformly to a function $\varphi$ which satisfies the same properties as $\varphi_i$.
	
	Let $p_i = (x_i, z_i) \in \RR^{d-1} \times \RR$. Then in particular
	\begin{equation}\label{tmp:diffvertdist}
		d_{p_i} \leq z_i - \varphi_i(x_i) \leq d_{p_i} \sqrt{ 1+\left|\theta\left(\frac32 R \right) \right|^2} + \theta(r_i) d_{p_i} \leq 2d_{p_i}.
	\end{equation} 
	Simple computations show that $\frac{D_i - p_i}{r_i} \cap B_5(0)$ corresponds to the region above the graph of the function $\psi_i: \RR^{d-1} \to \RR$, defined as
	\[ \psi_i(y) := \frac{1}{r_i} \left[ \varphi_i(x_i+r_i y) - \varphi_i(x_i) - (z_i - \varphi_i(x_i)) \right]. \]
	By \eqref{tmp:diffvertdist} and the assumption \eqref{assp:distr}, we have that $\psi_i(y)$, modulo passing to a subsequence, converges uniformly to a linear function $\varphi_\infty(y):= \langle \nabla \varphi(x_\infty), y \rangle$, where $x_\infty$ is a cluster point for $\{x_i\} \subset B_R(0)$. In other words, inside $B_5(0)$ the sequence of domains
	\[ \frac{D_i-p_i}{r_i} \text{ converges to the upper half space } D_\infty:= \{(y,y_d)\in \RR^{d-1} \times \RR: y_d > \varphi_\infty(y) \}.  \]
	
	By Lemma \ref{lm:urunifgrad}, $T_{p_i,r_i}u_i(0) = 0$
	and compactness, we get that the sequence 
	\[ T_{p_i, r_i} u_i(Y): = \frac{u_i(p_i + r_i Y)-u_i(p_i) }{ \left( \frac{1}{r_i^d} \mathlarger{\iint}_{B_{r_i}(p_i) \cap D_i} |u_i-u_i(p_i)|^2 ~dZ \right)^{\frac12} } \] 
	converges uniformly and in $W^{1,2}$ to a harmonic function $u_\infty$ in $D_\infty \cap B_5(0)$. 
	On the other hand, the sequence 
	\[ T_{q_i, r_i} u_i(Y) := \frac{u_i(q_i+r_i Y)}{\left(\frac{1}{r_i^d} \mathlarger{\iint}_{B_{r_i}(q_i)} u_i^2 ~dZ \right)^{\frac12} } 
	\]
	(and the sequence $T_{r_i} v_{q_i}$) also converges uniformly and in $W^{1,2}$ to a harmonic function $\tilde{u}_\infty$ in the same upper half space $D_\infty$.
	Moreover we claim there exists some constant $a_i$ which is uniformly bounded from above and below such that
	\begin{equation}\label{claim:changebin}
		|a_i T_{p_i,r_i}u_i(Y) - T_{q_i, r_i} u_i(Y)| \to 0 \text{ uniformly for } Y \in B_5(0).
	\end{equation}
	Assuming the claim is true, then $\tilde{u}_\infty = u_\infty$. (A priori $\tilde{u}_\infty$ is a constant multiple of $u_\infty$, and the constant must be $1$ since they both have unit $L^2$ norm on $B_1(0)$.)
	Hence
	\begin{align*}
		\left|N_{\mathcal{C}}(u_i, p_i, r_i) - \widetilde{N}(v_{q_i}, r_i) \right| & = \left| N_{\mathcal{C}}(T_{p_i, r_i} u_i, 0, 1)  - N(T_{r_i} v_{q_i}, 0, 1) \exp \left(C\int_0^{r_i} \frac{\theta(4s)}{s} ~ds \right) \right| \\
		& \to \left| N(u_\infty, 0, 1) - N(\tilde{u}_\infty, 0,1) \right| = 0,
	\end{align*}
	which contradicts the assumption \eqref{cl:freqbin}.
%
%
	
	\textit{Proof of the claim \eqref{claim:changebin}}.
	By Lemma \ref{lm:uup}, Remark \ref{rm:lmuup} and that $\dist(p_i, \partial D_i) \ll r_i$ (by \eqref{assp:distr}, the continuity of $\tilde{\theta}$, and $r_i \to 0$), we have
	\[ \frac{1}{r_i^d} \iint_{B_{r_i}(p_i)} |u_i - u_i(p_i)|^2 ~dZ \approx \frac{1}{r_i^{d-1}} \int_{\partial B_{r_i}(p_i)} |u_i -u_i(p_i)|^2 \sH. \]
	Since $|p_i-q_i| \leq 2d_{p_i} \ll r_i$, as in the proof of Remark \ref{rm:lmuup} we have
	\[ \frac{1}{r_i^{d-1}} \int_{\partial B_{r_i}(p_i)} |u_i -u_i(p_i)|^2 \sH \approx \frac{1}{r_i^d} \iint_{B_{r_i}(q_i)} u_i^2 ~dZ. \]
	Let us denote
	\[ a_i =  \left(\frac{ \frac{1}{r_i^d} \mathlarger{\iint}_{B_{r_i}(p_i)} |u_i-u_i(p_i)|^2 ~dZ}{\frac{1}{r_i^d} \mathlarger{\iint}_{B_{r_i}(q_i)} u_i^2 ~dZ } \right)^{\frac12}. \]
	We have just shown that $a_i$'s are uniformly bounded from above and below.
%
%
	
	It remains to prove
	\begin{equation}\label{tmp:bint}
		\frac{u_i(p_i+r_i Y) - u_i(p_i)- u_i(q_i + r_i Y)}{\left(\frac{1}{r_i^d} \mathlarger{\iint}_{B_{r_i}(q_i)} u_i^2 ~dZ \right)^{\frac12} } \to 0.
	\end{equation}	
%
%
%
	By the sub-harmonicity of $u_i^2$ we have
	\[ |u_i(p_i ) |^2 \leq \fiint_{B_s(p_i)} u_i^2 \sH, \quad \text{ for any } s>0. \]
	For each $u_i$, we take $s=d_{p_i}$. Using $|p_i - q_i| \leq 2d_{p_i}$, the doubling property and the assumption \eqref{assp:distr}, we get
	\begin{align*}
		\frac{|u_i(p_i)|^2}{\frac{1}{r_i^d} \mathlarger{\iint}_{B_{r_i}(q_i)} u_i^2 ~dZ }\lesssim \frac{\frac{1}{d^d_{p_i}} \mathlarger{\iint}_{B_{3d_{p_i} }(q_i)} u_i^2 ~dZ }{\frac{1}{r_i^d} \mathlarger{\iint}_{B_{r_i}(q_i)} u_i^2 ~dZ } \lesssim \left(\frac{3d_{p_i}}{r_i} \right)^{2N_{q_i} \cdot \exp\left(-C\int_0^{r_i} \frac{\theta(s)}{s} ~ds \right)} \leq \left( \frac{3d_{p_i}}{r_i}\right)^{3/2} \to 0,
	\end{align*}
	where we have used $N_{q_i} \geq 1$ proven in \cite[Lemma 5.23]{KZ}.
	Next, by the boundary gradient estimate in \cite[Theorem 1.4.3]{CK}, the monotonicity formula at $q_i$, and the assumption \eqref{assp:distr}, we get
	\begin{align}
		\frac{\left| u_i(p_i+r_i Y) - u_i(q_i + r_i Y) \right|}{\left(\frac{1}{r_i^d} \mathlarger{\iint}_{B_{r_i}(q_i)} u_i^2 ~dZ \right)^{\frac12} } & \leq \frac{\left| p_i - q_i \right| }{\left(\frac{1}{r_i^d} \mathlarger{\iint}_{B_{r_i}(q_i)} u_i^2 ~dZ \right)^{\frac12} } \sup_{B_{\frac{11}{2} r_i}(q_i) \cap D_i } |\nabla u_i| \nonumber \\
		&\leq \frac{\left| p_i - q_i \right| }{\left(\frac{1}{r_i^d} \mathlarger{\iint}_{B_{r_i}(q_i)} u_i^2 ~dZ \right)^{\frac12} } \cdot \frac{1}{r_i} \left( \frac{1}{r_i^d} \mathlarger{\iint}_{B_{6r_i}(q_i)} u_i^2 ~dZ \right)^{\frac12} \nonumber \\
		& \lesssim \frac{|p_i-q_i|}{r_i} \left(\frac{ \frac{1}{r_i^d} \mathlarger{\iint}_{B_{6r_i}(q_i)} u_i^2 ~dZ  }{\frac{1}{r_i^d} \mathlarger{\iint}_{B_{r_i}(q_i)} u_i^2 ~dZ } \right)^{\frac12} \nonumber \\
		& \lesssim \frac{d_{p_i}}{r_i} \cdot 6^{N_{q_i}(6r_i) \exp\left(-C\int_0^{r_i} \frac{\theta(s)}{s} ~ds \right)} \nonumber \\
		& \lesssim \tilde{\theta}\left(\frac{r_i}{\rho} \right) \cdot 6^{C(\Lambda)} \to 0,\label{tmp:bint2}
	\end{align}
	since $r_i \to 0$.
	This finishes the proof of \eqref{claim:changebin}.
\end{proof}

\section{Proof of Theorem \ref{thm:main}}

In Section \ref{sec:domtrsf}, in a neighborhood of any point $X_0 \in \pD \cap B_{\frac{R}{40}}(0)$ we defined $\tilde{u} = u \circ G$ on $\RR^d_+ \cap B$ and extended $\tilde{u}$ to the entire ball $B$ by an odd reflection; moreover, we have shown that $\tilde{u}$ is a solution to the elliptic equation $\divg(A(\cdot) \nabla \tilde{u}) = 0$ in $B$, where the coefficient matrix $A(\cdot)$ is $\alpha$-H\"older regular, with constant that only depends on the corresponding constant $C_\alpha$ for the $C^{1,\alpha}$ domain $D$. We claim that $\tilde{u}$ satisfies the normalized doubling assumption \eqref{cond:doubling} in $B$. In fact, since $\tilde{u}$ is defined by an odd reflection across $\RR^d_+$, it suffices to show that $\tilde{u}$ satisfies the normalized doubling assumption in the half ball $\RR^d_+ \cap B$. This follows immediately from the normalized doubling property of $u$ proven in Corollary \ref{cor:doublingu_in}, with a constant depending on $d, R$ and $ \Lambda$, by recalling that the domain transformation satisfies
\begin{equation}\label{eq:DGbd}
	\det DG \approx 1 \quad \text{ on } \RR^d_+ \cap B. 
\end{equation} 
Therefore, we may apply Theorem \ref{thm:HJ} to $\tilde{u}$ and conclude that 
\[ \mathcal{H}^{d-2}(\mathcal{C}(\tilde{u}) \cap B/2) \leq C, \]
where the constant $C$ depends on $d, R, \Lambda$ as well as $\alpha, C_\alpha$.
Again by \eqref{eq:DGbd}, we have that
\[ \mathcal{H}^{d-2}( \mathcal{C}(u) \cap G(\overline{\RR^d_+} \cap B/2)) \approx \mathcal{H}^{d-2}(\mathcal{C}(\tilde{u}) \cap \overline{\RR^d_+} \cap B/2) \leq C. \]

Recall that the radius of the ball $B$ only depends on the modulus of continuity of $\nabla \varphi$. More precisely when $D$ is a $C^{1,\alpha}$ domain with constant $C_\alpha$, the radius of $B$ only depends on $\alpha$ and $C_\alpha$.
 A simple covering argument says that we can find finitely many transformation maps $G_i$ and balls $B_i$'s, $i=1, \cdots, N$, so that
\[ \overline{D} \cap B_{\frac{R}{50}}(0) \subset \bigcup_{i=1}^N G_i(\overline{\RR^d_+} \cap B_i/2), \]
and the constant $N$ only depends on $\alpha, C_\alpha$ and $R$. Therefore we conclude that
\[ \mathcal{H}^{d-2}(\mathcal{C}(u) \cap B_{\frac{R}{50}}(0)) \leq \sum_{i=1}^N \mathcal{H}^{d-2}\left( \mathcal{C}(u) \cap G_i( \overline{\RR^d_+} \cap B_i/2 \right) \leq C'.  \]
This finishes the proof of Theorem \ref{thm:main}.

\appendix
\renewcommand{\theequation}{A.\arabic{equation}} 
\setcounter{equation}{0}  
\section*{Appendix}

It is tempting to apply the arguments in \cite{KZ} to the critical set instead of the singular set, by replacing the frequency function $N_{\mathcal{S}}(\cdot)$ by the \emph{normalized} frequency function $N_{\mathcal{C}}(\cdot)$, which is more suitable to study the critical set. However, the dominating terms in the derivatives of $N_{\mathcal{C}}(\cdot)$ and $N_{\mathcal{S}}(\cdot)$ are different (compare the expressions of $R_h(r)$ in \cite[(3.14)]{KZ} and \eqref{def:Rhr}); and the precise form of the derivative of the frequency function is used in a fundamental way to control the $(d-2)$-dimensional $\beta$-number by the change of the frequency function, see \cite[Section 9]{KZ}, which in turn allows us to use Reifenberg-type theorems to conclude the covering argument. 
More precisely, to include all critical points, we would still be able to estimate the $\beta$-number $\beta_{\mu}^{d-2}(p,r)$ from above by the changes of the \emph{normalized} frequency function as in \cite[Theorem 9.7]{KZ}, but we would get an extra term on the right hand side of \cite[(9.8)]{KZ}, which is of the form
\[ \Err(p,r) := \frac{\frac{1}{r^{d-2}} \mathlarger{\int}_{B_r(p)} |u(X)-u(p)|^2 ~d\mu(X) }{\frac{1}{r^{d}} \mathlarger{\iint}_{B_r(p)} |u-u(p)|^2 ~dZ }.  \]

In fact, it is not hard to show that $\Err(p,r)$ is small, since the oscillation of the values of $u$ on the critical set is indeed small (in a scale-invariant way); however, not only do we need $\Err(p,r)$ to be small, we also need it to be summable across the scale $r$, in order to apply Reifenberg-type theorems. We set out to prove this in Lemma \ref{lm:spvar} for harmonic functions in the interior; moreover, following the arguments laid out above, Lemma \ref{lm:spvar} gives an alternative way to estimate the size of the critical set for harmonic functions in the interior. (Unfortunately, when we consider its dependence on the upper bound of the frequency function $\Lambda$, this alternative proof only gives $\exp(C(\Lambda))$, which could be worse than the bound $C(d)^{\Lambda^2}$ proven in \cite[Theorem 1.1]{NVCS}. This is mainly because the constants in Lemma \ref{lm:spvar} depend on $\Lambda$.) In the following, unless specified otherwise, the notation $N(p,r) = N_{\mathcal{C}}(p,r)$ denotes the \emph{normalized} frequency function as defined in \eqref{def:DrHrin} and \eqref{def:Nrin}.

\begin{lemma}\label{lm:spvar}
	Let $u$ be a harmonic function in $B_3(0) \subset \RR^d$ such that $N(0,3) \leq \Lambda$. Then there exists a constant $C$ (depending on $d$ and $\Lambda$), such that for any $0<r<1$ and any $X_1, X_2\in B_1(0)$ with $|X_1 - X_2| \leq r/2$, we have
	\begin{equation}\label{eq:spvar}
		|N(X_1, r) - N(X_2, r)| \leq C \left( W^{\frac12}(X_1, r) + W^{\frac12}(X_2, r) \right),
	\end{equation}
	where $W(X_j, r) = N(X_j, 3r/2) - N(X_j, r/2)$ is assumed to be small, for $j=1,2$.
	
	Moreover, if additionally $X_1 \in \Ct(u)$, we also have 
	\begin{equation}\label{eq:svvalue}
		|u(X_1) - u(X_2)| \leq C\left( W^{\frac12}(X_1, r) + W^{\frac12}(X_2, r) \right) \cdot \left( \fiint_{B_r(X_1)} |u-u(X_1)|^2 ~dZ \right)^{\frac12}.
	\end{equation}
	The constants $C$ above depend on $d$ and $\Lambda$.
\end{lemma}
\begin{remark}
\begin{itemize}
	\item The estimate \eqref{eq:svvalue} is clearly not true if $X_1$ does not belong to the critical set (or the effective critical set, see the definition in \cite[(2.14)]{KZ}), by simply considering $u$ being a linear function in the $(X_2-X_1)$-direction.
	\item The first statement of the above lemma in \eqref{eq:spvar} is analogous to \cite[Lemma 8.1]{KZ}, the spatial variation of the usual frequency function.
\end{itemize}
	
\end{remark}
\begin{proof}
	Let $Y = X_2 - X_1$. Any point on the line segment $[X_1, X_2]$ can be written as $X_t=X_1 + t Y$, where $t\in [0,1]$. 
	
	\textbf{Step 1. Proof of \eqref{eq:spvar}.} Let 
	\[ \rho_{X_t}(Z) = \frac{Z-X_t}{|Z-X_t|} \]
	denote the radial direction with vertex at $X_t$, and let $\partial_{\rho_{X_t}}$ denote the derivative in the radial direction $\rho_{X_t}$. We have
	\begin{align*}
		\frac{d}{dt} H(X_t,r) = 2\int_{\partial B_r(X_t)} \langle \nabla u, Y \rangle (u-u(X_t)) ~d\mathcal{H}^{d-1}- 2 \langle \nabla u(X_t), Y \rangle \int_{\partial B_r(X_t) \cap D} (u-u(X_t)) ~d\mathcal{H}^{d-1},
	\end{align*}
	\begin{align*}
		\frac{d}{dt} D(X_t,r) = \int_{\partial B_r(X_t)} \langle \nabla u, Y \rangle \partial_{\rho_{X_t}} u \sH. 
	\end{align*}
	Hence
	\begin{align}
		& \frac{d}{dt} N(X_1 + tY, r) \nonumber \\
		& = \frac{d}{dt} \log N(X_1 + tY, r) \cdot N(X_1 + tY, r) \nonumber \\
		& = \frac{2}{H(X_t, r)} \left[ r \int_{\partial B_r(X_t)} \langle \nabla u, Y \rangle ~\partial_{\rho_{X_t}} u  ~d\mathcal{H}^{d-1} - N(X_t, r) \cdot \int_{\partial B_r(X_t)} \langle \nabla u, Y \rangle \left( u - u(X_t) \right) ~d\mathcal{H}^{d-1} \right] \nonumber \\
		& = \frac{2}{H(X_t, r)} \left[ \int_{\partial B_r(X_t)} \langle \nabla u, Y \rangle ~\left[ r \partial_{\rho_{X_t}} u  - N(X_t, r) \cdot \left( u - u(X_t) \right) \right] ~d\mathcal{H}^{d-1} \right].\label{tmp:spvarderNr}
	\end{align} 
	Let 
	\[ \E_j(Z) = \langle \nabla u(Z), Z-X_j \rangle - N(X_j, |Z-X_j|) \cdot \left( u(Z)-u(X_j) \right), \quad \text{ for } j=1,2.  \]
	Then
	\begin{align}
		& \langle \nabla u(Z), Y\rangle \nonumber \\
		& = \E_1(Z) - \E_2(Z) + N(X_1, |Z-X_1|) \cdot \left( u(Z)-u(X_1) \right) - N(X_2, |Z-X_2|) \cdot \left( u(Z)-u(X_2) \right) \nonumber \\
		& = \E_1(Z) - \E_2(Z) + N(X_1, r) \cdot \left( u(Z)-u(X_1) \right) - N(X_2, r) \cdot \left( u(Z)-u(X_2) \right) \nonumber \\
		& \quad + \left[N(X_1, |Z-X_1|) - N(X_1, r) \right]\left( u(Z)-u(X_1) \right) - \left[N(X_2, |Z-X_2|) - N(X_2, r)\right]\left( u(Z)-u(X_2) \right) \nonumber \\
		& = \left[ \E_1(Z) - \E_2(Z)\right] + \left[N(X_1, r) - N(X_2, r) \right] \cdot \left( u(Z)-u(X_t) \right) \nonumber \\
		& \quad + \left[ N(X_1, r) \cdot \left(u(X_t) - u(X_1) \right) -N(X_2, r) \cdot \left(u(X_t) - u(X_2) \right) \right] \nonumber \\
		& \quad + \left[N(X_1, |Z-X_1|) - N(X_1, r) \right]\left( u(Z)-u(X_1) \right) - \left[N(X_2, |Z-X_2|) - N(X_2, r) \right]\left( u(Z)-u(X_2) \right).\label{tmp:gradu}
	\end{align}
	
	To calculate the second and the third terms, notice that 
	\begin{align*}
		&  \int_{\partial B_r(X_t)} r \left( u-u(X_t) \right) \partial_{\rho_{X_t}} u ~d\mathcal{H}^{d-1} - N(X_t,r) \int_{\partial B_r(X_t)} \left( u-u(X_t) \right)^2 ~d\mathcal{H}^{d-1}\\
		& \qquad \qquad =  r D(X_t, r) 
		- N(X_t,r) H(X_t, r) = 0,
	\end{align*}
	and 
		\begin{align*}
		& \int_{\partial B_r(X_t)} \left[ r \partial_{\rho_{X_t}} u  - N(X_t, r) \cdot \left( u - u(X_t) \right) \right] ~d\mathcal{H}^{d-1} \\
		& = r\iint_{B_r(X_t)} \divg(\nabla u) ~ dZ 
		- N(X_t, r) \cdot \int_{\partial B_r(X_t)} \left( u-u(X_t) \right) ~d\mathcal{H}^{d-1} = 0,
	\end{align*}
	where we use the mean value property in the last equality. 
	Inserting these equalities into \eqref{tmp:gradu} and \eqref{tmp:spvarderNr} we get
	\begin{align}
		& \frac{d}{dt} N(X_1 + tY, r) \nonumber \\
		& = \frac{2}{H(X_t, r)} \int_{\partial B_r(X_t)} \left( \E_1- \E_2 \right) \left[ r\partial_{\rho_{X_t}} u - N(X_t,r) \cdot \left( u-u(X_t) \right) \right] d\mathcal{H}^{d-1} \nonumber \\
		& + \frac{2}{H(X_t,r)} \int_{\partial B_r(X_t)} \left( N(X_1, |Z-X_1|) - N(X_1, r) \right) (u - u(X_1))\left[ r\partial_{\rho_X} u - N(X,r) \cdot \left( u-u(X_t) \right) \right] d\mathcal{H}^{d-1} \nonumber \\
		& - \frac{2}{H(X_t,r)} \int_{\partial B_r(X_t)} \left( N(X_2, |Z-X_2|) - N(X_2, r) \right) (u - u(X_2))\left[ r\partial_{\rho_X} u - N(X,r) \cdot \left( u-u(X_t) \right) \right] d\mathcal{H}^{d-1} \nonumber \\
		& = : \I_0(X_t) + \I_1(X_t)  + \I_2(X_t) \label{tmp:spvar}
	\end{align}

	For any $X_t$ on the line segment $[X_1, X_2]$ and any $Z\in \partial B_r(X_t)$, we have
	\[ \frac{r}{2} \leq |Z-X_j| \leq \frac32 r, \quad \text{ for } j=1,2. \]
	Hence
	\[ 
	\left| N(X_j, |Z-X_j|) - N(X_j, r) \right| \leq W(X_j, r), \]
	by the monotonicity of $N(X_j, \cdot)$.
	We can estimate the remaining terms of $\I_j(X_t), j=1, 2$ as follows:
	\begin{align}
		| \I_j(X_t) | & \leq \frac{2}{H(X_t,r)}  W(X_j, r) \int_{\partial B_r(X_t)} \left|(u - u(X_j))\left[ r\partial_{\rho_{X_t}} u - N(X_t,r) \cdot \left( u-u(X_t) \right) \right] \right| d\mathcal{H}^{d-1} \nonumber \\
		& \leq \frac{1}{H(X_t,r)}  W(X_j, r)\left[ \left( r + 2N(X_t, r)\right) \cdot \int_{\partial B_r(X_t)} \left(u-u(X_j) \right)^2 ~d\mathcal{H}^{d-1}  \right. \nonumber \\
		& \qquad \left. + r \int_{\partial B_r(X_t)} |\nabla u|^2 ~d\mathcal{H}^{d-1} \right] \nonumber \\
		& \leq   W(X_j, r) \left[ C(\Lambda) \frac{\int_{\partial B_r(X_t)} \left(u-u(X_j) \right)^2 ~d\mathcal{H}^{d-1} }{H(X_t,r) } + \frac{r\int_{\partial B_r(X_t)} |\nabla u|^2 ~\mathcal{H}^{d-1} }{H(X_t,r) } \right]. \label{tmp:I2}
	\end{align}
	We estimate $\I_0(X)$ by the Cauchy-Schwarz inequality, and obtain
	\begin{equation}\label{tmp:I1}
		|\I_0(X_t)| \lesssim \left[ r \left( \frac{ \int_{\partial B_r(X_t)} |\nabla u|^2 ~d\mathcal{H}^{d-1}}{H(X_t,r)} \right)^{\frac12} + N(X_t,r) \right]\cdot  \left( \frac{1}{H(X_t,r)} \int_{\partial B_r(X_t)} |\E_1(Z)|^2 + |\E_2(Z)|^2 ~d\mathcal{H}^{d-1} \right)^{\frac12}.
	\end{equation}
	
	To further estimate $\I_0(X_t)$ and $\I_j(X_t)$, we make the following observation. Since $u$ is a harmonic function, $(u-c)^2$ is sub-harmonic for any constant $c$. Hence for every $X$,
	\[ \fint_{\partial B_\rho(X)} \left(u-u(X) \right)^2 \sH \text{ is increasing with respect to }\rho. \]
	We also have that 
	\begin{align*}
		\left|u(X_j) -u(X_t) \right| = \left|\fiint_{B_r(X_j)} \left( u-u(X_t) \right) ~dZ \right| & \leq \left( \fiint_{B_r(X_j)} \left( u-u(X_t) \right)^2 ~dZ\right)^{\frac12} \\	
		& \lesssim \left( \fiint_{B_{2r}(X_t)} \left( u-u(X_t) \right)^2 ~dZ\right)^{\frac12}.
	\end{align*}
	For any $r_1<r_2$, let $A_{r_1, r_2}(X) := B_{r_2}(X) \setminus B_{r_1}(X)$ denote the annulus centered at $X$.
	It follows that 
	\begin{align*}
		H(X_j, r) = \int_{\partial B_r(X_j)} \left(u-u(X_j) \right)^2 \sH & \leq \frac2r \iint_{A_{r, \frac32 r}(X_j)} \left(u-u(X_j) \right)^2 ~dZ \\
		& \leq \frac2r \iint_{B_{2r}(X_t)} \left(u-u(X_j) \right)^2 ~dZ \\
		& \lesssim \frac1r \iint_{B_{2r}(X_t)} \left(u-u(X_t) \right)^2 ~dZ \\
		& \lesssim \int_{\partial B_{2r}(X_t)} \left(u-u(X_t) \right)^2 \sH = H(X_t, 2r).
	\end{align*}
	Combined with the doubling property, we have that 
	\begin{equation}\label{eq:Hchangecenter}
		H(X_j, r) \lesssim H(X_t,2r) \lesssim_{\Lambda} H(X_t,r), \quad \text{ for }j=1,2.
	\end{equation} 
	Therefore modulo constants, we can change the denominators in \eqref{tmp:I2} and \eqref{tmp:I1} to $H(X_j, r)$.
	
	Now we allow $X=X_t$ to move in the line segment $[X_1, X_2]$, as  $t$ varies in $[0,1]$. Since $|X_1 - X_2| \leq r/2$, the integration region satisfies
	\[ \bigcup_{t\in [0,1]} \partial B_r(X_t) \subset A_{\frac{r}{2}, \frac32 r}(X_j). \] 
	By \eqref{tmp:I2} and \eqref{eq:Hchangecenter}, we have
	\begin{align}
		\int_0^1 |\I_j(X_t)| ~dt & \leq W(X_j, r) \left[  \frac{C_1}{H(X_j,r) } \iint_{A_{\frac{r}{2}, \frac32 r} (X_j) } \left( u-u(X_j) \right)^2 ~dZ + \frac{C_2 r}{H(X_j, r)} \iint_{A_{\frac{r}{2}, \frac32 r} (X_j) } |\nabla u|^2 ~dZ \right] \nonumber \\
		& \lesssim W(X_j,r) \left[ C'_1+  C'_2 \frac{r D(X_j, \frac32 r)}{H(X_j, \frac32 r)}  \right] \nonumber \\
		& \leq C(\Lambda)W(X_j, r).\label{tmp:I2f}
	\end{align}
	Next we use the Cauchy-Schwarz inequality to estimate
	\begin{align*}
		\int_0^1 |\I_0(X_t)| ~dt &\lesssim \left\{ \int_0^1 \left[ \frac{r^2}{H(X_t, r)} \int_{\partial B_r(X_t)} |\nabla u|^2 \sH + N^2(X_t,r) \right] ~dt \right\}^{\frac12}  \\
		& \qquad \times \left\{ \int_0^1 \frac{1}{H(X_t, r)} \int_{\partial B_r(X_t)} |\E_1(Z)|^2 + |\E_2(Z)|^2 \sH ~dt \right\}^{\frac12} \\
		& \lesssim \left[ \frac{r^2}{H(X_j, \frac32 r)} \cdot  D(X_j, \frac32 r) + C(\Lambda) \right]^{\frac12} \\
		& \qquad \times \left[ \sum_{j=1}^2 \frac{1}{H(X_j,  \frac32 r)} \iint_{A_{\frac{r}{2}, \frac32 r}(X_j)} |\E_j(Z)|^2 ~dZ \right]^{\frac12}.
	\end{align*}
	Recall Proposition \ref{prop:intmonotonicity} (the interior case, i.e. when $r< \dist(p,\pD)$), in particular \eqref{eq:dernNr} and \eqref{def:Rhr}. We have
	\begin{align}
		& \frac{1}{H(X_j, \frac32 r)} \iint_{A_{\frac{r}{2}, \frac32 r}(X_j)} |\E_j(Z)|^2 ~dZ \nonumber \\
		& = \frac{1}{H(X_j, \frac32 r)} \iint_{A_{\frac{r}{2}, \frac32 r}(X_j)} \left|\langle \nabla u(Z), Z-X_j \rangle - N(X_j, |Z-X_j|) \cdot \left(u(Z)-u(X_j) \right) \right|^2 ~dZ \nonumber \\
		& \lesssim r \int_{\frac{r}{2}}^{\frac32 r} \frac{\rho}{H(X_j, \rho)} \int_{\partial B_\rho(X_j)} \left|\partial_{\rho_{X_j}} u(Z) - \frac{N(X_j, \rho)}{\rho} \cdot \left(u(Z) - u(X_j) \right)  \right|^2 \sH ~d\rho \nonumber \\
		& \lesssim r \int_{\frac{r}{2}}^{\frac32 r} R_h(X_j, \rho) ~d\rho \nonumber \\
		& \lesssim W(X_j, r).\label{eq:estEbyW}
	\end{align}
	Therefore
	\begin{equation}\label{tmp:I1f}
		\int_0^1 |\I_0(X_t)| ~dt \leq C(\Lambda) \left( W^{\frac12}(X_1, r) + W^{\frac12}(X_2, r) \right). 
	\end{equation} 
	
	Combining \eqref{tmp:spvar}, \eqref{tmp:I2f} and \eqref{tmp:I1f}, we get
	\begin{equation}\label{eq:svfreq}
		|N(X_2, r) - N(X_1, r)| \leq \int_0^1 |\I_0(X_t)| + |\I_1(X_t)| + |\I_2(X_t)| ~dt \leq C_{sv} \left(W^{\frac12}(X_1, r) + W^{\frac12}(X_2, r) \right),
	\end{equation}
	where the constant $C_{sv}$ depends on $d$ and $\Lambda$.
	
	\textbf{Step 2. Proof of \eqref{eq:svvalue}.}
	Using the same notation, we have that
	\begin{align*}
		u(X_2) - u(X_1) =  \int_0^1 \frac{d}{dt} u(X_t) ~dt  = \int_0^1 \langle \nabla u(X_t), Y \rangle ~dt = \int_0^1 \fint_{\partial B_r(X_t)} \langle \nabla u, Y \rangle \sH ~dt ,
	\end{align*}
	As above, we write
	\begin{align*}
		& \fint_{\partial B_r(X_t)} \langle \nabla u, Y \rangle \sH \\
		& = \fint_{\partial B_r(X_t)} \E_1(Z) - \E_2(Z) + N(X_1, |Z-X_1|) \cdot \left( u(Z) - u(X_1) \right) -N(X_2, |Z-X_2|) \cdot \left( u(Z) - u(X_2) \right) \sH \\
		& = \fint_{\partial B_r(X_t)} \E_1(Z) - \E_2(Z) \sH + \fint_{\partial B_r(X_t)} \bar N \cdot (u(Z)-u(X_1)) - \bar N \cdot (u(Z) - u(X_2)) \sH \\
		& \qquad + \fint_{\partial B_r(X_t)} \left[ N(X_1, |Z-X_1|) - \bar{N} \right]\cdot (u(Z) - u(X_1)) \sH \\
		& \qquad - \fint_{\partial B_r(X_t)} \left[ N(X_2, |Z-X_2|) - \bar{N} \right]\cdot (u(Z) - u(X_2)) \sH \\
		& = \fint_{\partial B_r(X_t)} \E_1(Z) - \E_2(Z) \sH + \bar N \cdot \left( u(X_2) - u(X_1) \right) \\
		& \qquad + \fint_{\partial B_r(X_t)} \left[ N(X_1, |Z-X_1|) - \bar{N} \right]\cdot (u(Z) - u(X_1)) \sH \\
		& \qquad - \fint_{\partial B_r(X_t)} \left[ N(X_2, |Z-X_2|) - \bar{N} \right]\cdot (u(Z) - u(X_2)) \sH,
	\end{align*}
	where $\bar{N}$ is a constant to be chosen.
	Hence moving the term $\bar{N} \cdot (u(X_2) - u(X_1))$ to the left hand side, we can estimate it as follows:
	\begin{align*}
		& \left( \bar{N} - 1 \right) \cdot |u(X_2) - u(X_1)| \leq \frac{1}{|\partial B_r|} \left[ \int_0^1 \int_{\partial B_r(X_t)} |\E_1 - \E_1 | \sH ~dt \right. \\
		& \qquad \left. + \sum_{j=1}^2 \sup_{\partial B_r(X_t)} \left|N(X_j, |Z-X_j|) - \bar{N} \right| \int_0^1 \int_{\partial B_r(X_t)} |u-u(X_j)| \sH ~dt \right].
	\end{align*}
	Let $\bar{N} = N(X_1, r)$. Then by \eqref{eq:svfreq} we have
	\[ |N(X_2, r) - \bar{N}| \leq C_{sv} \left(W^{\frac12}(X_1, r) + W^{\frac12}(X_2, r) \right). \]
	Thus
	\begin{align*}
		\sup_{\partial B_r(X_t)} \left|N(X_j, |Z-X_j|) - \bar{N} \right| & \leq \max\left\{W(X_1, r), W(X_2, r) + C_{sv} \left(W^{\frac12}(X_1, r) + W^{\frac12}(X_2, r) \right) \right\} \\
		& \lesssim W^{\frac12}(X_1, r) + W^{\frac12}(X_2, r).
	\end{align*} 
	On the other hand, since $X_1 \in \Ct(u)$, we have that $T_{X_1, \rho} u$ has sequential limit to a homogeneous harmonic polynomial of degree at least $2$. Thus
	\[ \bar{N} = N(X_1, r) \geq \lim_{\rho \to 0} N(X_1, \rho) \geq 2. \]
	Therefore
	\begin{align*}
		& |u(X_2) - u(X_1)| \\
		& \lesssim \frac{1}{|\partial B_r|} \left[ \sum_{j=1}^2 \iint_{A_{\frac{r}{2}, \frac32 r}(X_j)} |\E_j| ~dZ + \left(W^{\frac12}(X_1, r) + W^{\frac12}(X_2, r) \right) \sum_{j=1}^2 \iint_{A_{\frac{r}{2}, \frac32 r}(X_j)} |u-u(X_j)| ~dZ \right] \\
		& \lesssim r \left[ \sum_{j=1}^2 \left( \fiint_{A_{\frac{r}{2}, \frac32 r}(X_j)} |\E_j|^2 ~dZ \right)^{\frac12} + \left(W^{\frac12}(X_1, r) + W^{\frac12}(X_2, r) \right) \sum_{j=1}^2 \left( \fiint_{A_{\frac{r}{2}, \frac32 r}(X_j)} |u-u(X_j)|^2 ~dZ\right)^{\frac12} \right] \\
		& \lesssim  \left(W^{\frac12}(X_1, r) + W^{\frac12}(X_2, r) \right) \cdot \sum_{j=1}^2 \left( \fiint_{B_r(X_j)} |u-u(X_j)|^2 ~dZ\right)^{\frac12} \\
		& \lesssim \left(W^{\frac12}(X_1, r) + W^{\frac12}(X_2, r) \right) \left( \fiint_{B_r(X_1)} |u-u(X_1)|^2 ~dZ\right)^{\frac12},
	\end{align*}
	for $j=1,2$, where we used arguments as in \eqref{eq:estEbyW} in the last inequality.
\end{proof}

\end{document}